# Superposition of interacting stochastic processes with memory and its application to migrating fish counts[1]

Hidekazu Yoshioka [a, *]

[a] Japan Advanced Institute of Science and Technology, 1-1 Asahidai, Nomi, Ishikawa 923-1292, Japan
[*] Corresponding author: yoshih@jaist.ac.jp, ORCID: 0000-0002-5293-3246

**Abstract**
Stochastic processes with long memories, known as long memory processes, are ubiquitous in various science and engineering problems. Superposing Markovian stochastic processes generates a non-Markovian long memory process serving as powerful tools in several research fields, including physics, mathematical economics, and environmental engineering. We formulate two novel mathematical models of long memory process based on a superposition of interacting processes driven by jumps. The mutual excitation among the processes to be superposed is assumed to be of the mean field or aggregation form, where the former yields a more analytically tractable model. The statistics of the proposed long memory processes are investigated using their moment-generating function, autocorrelation, and associated generalized Riccati equations. Finally, the proposed models are applied to time series data of migrating fish counts at river observation points. The results of this study suggest that an exponential memory or a long memory model is insufficient; however, a unified method that can cover both is necessary to analyze fish migration, and our model is exactly the case.

**Funding statement:** This study was supported by the Japan Society for the Promotion of Science (KAKENHI Nos. 22K14441 and 22H02456) and the Japan Science and Technology Agency (PRESTRO, No. JPMJPR24KE).

**Data availability statement:** The corresponding author will make data available upon reasonable request.

**Competing interests:** The author has no competing interests.

**Declaration of generative AI in scientific writing:** The authors did not use generative AI technology to write the manuscript.

**Permission to reproduce material from other sources:** N.A.

**Author contributions:** The first author conducted all parts of this research.

---

[1] See also the correction note attached to this preprint that corrects **Proposition 3**.

**Acknowledgments:** The author thanks Dr. Yumi Yoshioka for kindly preparing **Figure 1**. The manuscript has been significantly improved with the help of careful reviews by anonymous reviewers.

# 1. Introduction

## 1.1 Study background

Stochastic process models are widely used in many research fields to describe randomly fluctuating phenomena (Chapters 4-7 of Capasso and Bakstein [1]). Autocorrelation is a statistic that characterizes a stochastic process and measures memory decay. Markovian processes, such as solutions to stochastic differential equations (SDEs) driven by white-noise processes [2,3], have exponential memory that decays rapidly over time. In contrast, many real-world phenomena involve longer memories that decay only polynomially—the long memory process (Chapter 1 of Beran et al. [4]). Long memory processes exist in air traffic [5], earthquake-induced geochemical variations [6], realized volatility in finance [7], air temperature trends due to global warming [8], epidemic dynamics [9], and superstatistics [10]. Several methods exist for generating long memory processes, including the fractional (or Volterra) and superposition approaches, which are understood through infinite-dimensional stochastic processes, as reviewed below. This study focuses on superposition but reviews both to help us better understand their similarities and differences.

The fractional approach assumes that the memory of a system, which is often identified with a fractional derivative in time, has a polynomial decay. Therefore, the long memory nature of the fractional approach originates from temporal fractional differentiation. The fractional Brownian motion is the simplest model in the fractional approach, as a long memory generalization of classical Brownian motion [11]. Advanced examples include the fractional Hawkes processes [12-14], particle transport in water flows [15], geophysics and biology [16], fractional multi-state models [17], and coupled SDEs in finance and economics [18,19].

By contrast, the superposition approach assumes that a long memory process arises from superposing (i.e., integrating) suitable mutually independent Markovian stochastic processes with respect to their reversion speeds. The reversion speed of a process is interpreted as the reciprocal of its timescale of variations. Therefore, the long memory of the superposition approach originates from the coexistence of slow to fast timescale processes. The simplest model in the superposition approach is the superposition of the Ornstein–Uhlenbeck processes (supOU process) [20]. Advanced examples include the superposed economic models for stock and asset dynamics [21], mass estimate of galactic black holes [22], random coefficient autoregressive models along with their superpositions [23], anomalous diffusion [24], noise-driven Hamiltonian systems [25], and time-periodic river flows [26]. Superstatistics that deal with a doubly random variable [27,28] are closely related to the superposition approach because a superposed process is based on a system of stochastic processes with randomized coefficients.

supOU process is a non-Markovian process, as polynomial decaying autocorrelation suggests [29]. Furthermore, Volterra processes are superposed processes driven by fully correlated and common noise [30-32]. From this perspective, the fractional and superposition approaches generate a long memory process by superposing infinite processes, and are deeply related with each other. A significant difference between them is the assumption regarding the processes to be superposed. The fractional approach superposes processes driven by fully correlated noise (e.g., univariate martingale), whereas the superposition approach uses independent processes driven by independent noises. The independent driving noise in the superposition approach allows for the explicit derivation of statistics (e.g., autocorrelation and moments) of long memory processes [33,34]. This remarkable property of the superposition approach was exploited in the present study.

As explained above, the classical superposition approach assumes that the superposed processes are independent. An expected question is whether stochastic processes driven by mutually independent noises can be superposed without losing analytical tractability. This issue has not yet been addressed to the best of the author's knowledge. A class of mutually interacting processes was proposed for the fractional approach to generate a Hawkes-type process with long memory [35]. Mutual interactions among processes would become relevant for applied modeling studies, such as collective behavior in biology [36,37], finance [38], and chemical reactions [39], which are also interesting topics from a mathematical perspective.

### 1.2 Aim and contribution

The aims of this study, which address the limitations mentioned previously, are two-fold.

- ✓ Modeling long memory processes using mutually interacting processes based on the superposition approach.
- ✓ Application of the proposed models to biological count data.

The following paragraphs explain the contributions made toward achieving these aims. The focus is on processes driven by positive jumps, where mutual interaction is accounted for in state-dependent jump rates, different from previous models. More specifically, two types of mutual interactions are considered: mean-field and aggregated interactions.

Concerning the first aim, the mean-field interaction-based (MF) model assumes that the jump rates of the processes to be superposed are linked through their statistical average. In this case, the superposition is an infinite-dimensional version of mean-field SDEs [40,41]; however, the proposed model is more analytically tractable. In particular, the autocorrelation and memory of the resulting long memory process are explicitly obtained when the influence of the mean-field effect is visible. This explicit nature facilitates applications because a model can be identified by fitting it to data using matching statistics.

In contrast, the aggregation-based (AG) model assumes that the jump rates of the processes to be superposed are linked without assuming a mean-field (or ergodic) ansatz. In this case, the superposition is an infinite-dimensional version of the multivariate jump processes [42,43]. The absence of the mean-field ansatz results in a lack of explicit formulae for high-order statistics, such as variance; however, the moment-generating function can be found by solving a generalized Riccati equation given by a partial integro-differential equation. The well-posedness of generalized Riccati equations has been proven by its structure, including the quasi-monotonicity of the coefficients, in the literature [44,45]. We deal with the generalized Riccati equation in a unified way along with its associated linearized versions where weak (mild) solutions are defined in a Banach space of integrable functions. We numerically compute statistics of the AG model by using these equations.

Regarding the second aim, the proposed models are applied to the biological count data of the migrating fish *Plecoglossus altivelis altivelis* (*P. altivelis*) in Japan. This diadromous fish species is distributed in Northeast Asia along the Japan and East China Seas. It has a one-year life cycle and exhibits upstream migration from the sea, reservoir, or lake (i.e., a large water body) to a connected river during spring [46]. The downstream larvae migration during autumn [47], biological growth and spawning [48], swimming speed in water currents [49], and harvest management [50] of this fish species have been extensively studied. However, the upstream migration, particularly the memory structure in migrating population counts, has not been well studied, except for that based on the supOU process [51]. Several studies have discussed that fish migration is not passive but collective behavior due to social cues [52-54]. Such phenomena arise from nonlinear interactions among the fish [55,56]. However, resolving these mechanisms becomes inefficient when focusing on macroscopic migration. Our model's framework is conceptual and restrictive because it only deals with the migration fish count, not the environmental cues that could affect migration. Nevertheless, a stochastic process model would be effective for such time series analysis. Additionally, the application in this study serves as a building block for modeling fish migration based on long memory processes that consider social cues.

The application of the model in this study considers five rivers in Japan and reveals that the parameters of the MF model can be determined along with their realizability conditions. Applying the MF model to biological count data is a novel contribution, and its comparison with the AG model is another contribution of this study. The aim of the application is the characterization of fish migration from a stochastic modeling viewpoint, not simulation. The AG model is also computationally applied to data using Monte Carlo simulations and a generalized Riccati equation. Consequently, this study contributes to the modeling and application of novel long memory processes.

The rest of this paper is organized as follows. **Section 2** reviews earlier models that serve as building blocks of the proposed models. **Section 3** describes the formulation and analysis of the proposed models. These models are applied to the biological count data, as described in **Section 4**. **Section 5** summarizes the study and presents the perspectives of our research. **Appendices** provide proofs. **The supplementary material** includes auxiliary data.

2. **Superposition approach**

Superposition-based long memory processes are presented herein. Each process is an infinite-dimensional limit of a finite-dimensional model, with the latter being easier to comprehend. The finite-dimensional model is presented in **Section A1**. The time $t$ is a real parameter, and we work with a complete probability space $(\Omega, \mathbb{F}, \mathbb{P})$ ($\Omega$: collection of all events, $\mathbb{F}$: filtration, and $\mathbb{P}$: probability function) (e.g., Chapter 1.1 in Øksendal and Sulem [57]; Section 3 in Gomez et al. [25]). The Poisson random measures are double-sided to manage stationary processes (e.g., Section 2 in Barndorff-Nielsen and Stelzer [33]). Each model is adapted to a filtration with respect to its associated Poisson random measure. $\mathbb{E}$ and $\mathbb{V}$ represent the expectation and variance, respectively. Autocorrelation with the time-lag $\tau \geq 0$ is represented by $\rho(\tau)$. The left limit of the stochastic process $X = (X_t)_{t \in \mathbb{R}}$ at time $t$ is expressed as $X_{t-}$.

### 2.1 Previous model without mutual interaction

#### 2.1.1 Nominal process

First, a nominal process $X = (X_t)_{t \in \mathbb{R}}$ with the reversion speed $r > 0$ satisfies the SDE

$$dX_t = (b - rX_t)dt + \int_{z=0}^{z=+\infty} \int_{u=0}^{u=rX_{t-}} zN(dt, dz, du), \ t \in \mathbb{R}, \quad (1)$$

with the source rate as $b > 0$, and $N$ is a Poisson random measure of $\mathbb{R} \times (0, +\infty) \times (0, +\infty)$ with the compensator $dt v(dz) du$. Additionally, $v$ is a Lévy measure with the density $v(dz) = \sigma(z) dz$ and a measurable function $\sigma : (0, +\infty) \to (0, +\infty)$ such that $M_1 \in (0,1) : M_m = \int_0^{+\infty} z^m \sigma(z) dz$ ($m = 0, 1, 2, ...$). The SDE (1), subject to a non-zero initial condition, admits a pathwise unique solution that is non-negative and right-continuous with left limits (Proposition 5.6 of Fu and Li [58] with $b(x) = b - rx$ and $h_1(x) = rx$, using the notations in this literature).

At a stationary state, using (1) for $t \in \mathbb{R}$, the following is obtained:

$$X_t = \frac{b}{r} + \int_{s=-\infty}^{s=t} \int_{z=0}^{z=+\infty} \int_{u=0}^{u=rX_{s-}} e^{-r(t-s)} zN(ds, dz, du). \quad (2)$$

The quantity $M_1$ represents the degree of self-excitedness as in the intensity of the Hawkes process [59], where the assumption $M_1 \in (0,1)$ in our context implies that the moment of jumps should not be large to guarantee $X$ stationarity. An elementary calculation obtains the stationary average as $\mathbb{E}[X_t] = \frac{b}{r(1-M_1)}$. Higher-order moments are also obtained analytically, e.g., $\mathbb{V}[X_t] = \frac{M_2 b}{2r(1-M_1)^2}$. The autocorrelation of the nominal process is $\rho(\tau) = e^{-r(1-M_1)\tau}$, revealing an exponential decay with the exponent proportional to $r$ and $1 - M_1$. Therefore, the autocorrelation decays more slowly for a more strongly self-exciting case.

The nominal process with the couple $(r, b)$ of reversion speed and source rate is formally represented as $X(r, b)$.

#### 2.1.2 Previous model

This subsection is based on the study by Yoshioka [34] but with an explanation using measure-valued notations, as in that of Gomez et al. [25], where the superposition mechanism becomes more visible. The superposed process, as a building block of the models proposed in this study, is a superposition of a continuum of independent nominal processes parameterized by the reversion speed $r$. We assume that the reversion speed $r$ follows a probability measure $\pi$ satisfying $R = \int_0^{+\infty} r^{-1}\pi(\mathrm{d}r) < +\infty$, indicating that $\pi$ is not too much singular at $r=0$.

The superposed process $Z = (Z_t)_{t \in \mathbb{R}}$ in the previous model [34] is given for $t \in \mathbb{R}$ as

$$Z_t = \int_{r=0}^{r=+\infty} X_t(r, b\pi(\mathrm{d}r)), \qquad (3)$$

which can be rewritten as

$$Z_t = \int_{r=0}^{r=+\infty} \left\{ \frac{b\pi(\mathrm{d}r)}{r} + \int_{s=-\infty}^{s=t} \int_{z=0}^{z=+\infty} e^{-r(t-s)} z N_{\mathrm{Pr}}(\mathrm{d}s, \mathrm{d}z, \mathrm{d}r) \right\}. \qquad (4)$$

Here, $N_{\mathrm{Pr}}$ is a Poisson random measure with the jump rate $rX_{t-}(r, b\pi(\mathrm{d}r))\nu(\mathrm{d}z)\mathrm{d}t$. The origin of $X_{t-}(r, b\pi(\mathrm{d}r))$ in the jump rate is uncovered by considering a finite-dimensional version as a weak limit of (4) (**Section A1**). Each process $X(r, b\pi(\mathrm{d}r))$ is understood as a measure-valued process parameterized by $r > 0$, implying that it is infinitesimally small to certain extent, while its integration (i.e., the right side of (3)) is not. The average and variance of the process $X(r, b\pi(\mathrm{d}r))$ are proportional to $\pi(\mathrm{d}r)$ and are at the order of $\mathrm{d}r$ if $\pi$ has a density. We have (see **Section A2**)

$$\mathbb{E}[Z_t] = \frac{b}{1-M_1}R \text{ and } \mathbb{V}[Z_t] = \frac{M_2 b}{2(1-M_1)^2}R. \qquad (5)$$

Particularly, the positivity of the variance $\mathbb{V}[Z_t]$ suggests that the process $Z$ is stochastic. Moreover, the autocovariance $\mathbb{C}[Z_{t+\tau}, Z_t]$ of $Z$ with lag $\tau \geq 0$ is given by

$$\mathbb{C}[Z_{t+\tau}, Z_t] = \frac{\mathbb{V}[Z_t]}{R} \int_0^{+\infty} \frac{1}{r} e^{-r(1-M_1)\tau} \pi(\mathrm{d}r). \qquad (6)$$

Then, the autocorrelation is obtained as

$$\rho(\tau) = \frac{\mathbb{C}[Z_{t+\tau}, Z_t]}{\mathbb{V}[Z_t]} = \frac{1}{R} \int_0^{+\infty} \frac{1}{r} e^{-r(1-M_1)\tau} \pi(\mathrm{d}r). \qquad (7)$$

The previous model admits exponentially- (i.e., $\pi$ is a Dirac Delta) and polynomially decaying autocorrelations (i.e., $\pi$ is a Gamma distribution, as assumed in **Section 4**). Moreover, if $\nu$ is proportional to a probability measure, then $M_0 = \int_0^{+\infty} \nu(\mathrm{d}z) < +\infty$, and the total number of jumps $\mathbb{J}[Z_t]$ of $Z$ at a given unit time interval is given by

$$\mathbb{J}[Z_t] = \frac{M_0}{1-M_1}b. \qquad (8)$$

Model flexibility and the explicit availability of statistics are advantages of superposition.

**Remark 1:** In the remainder of this paper, equations and inequalities for $r > 0$ are understood almost surely (a.s.) with respect to the probability measure $\pi$.

**Remark 2:** The superposed process $Z$ is non-Markovian because it cannot be predicted based on the current value. Predicting $Z$ at a future time needs the information of $X(r,b\pi(\mathrm{d}r))$ for all $r>0$. Qualitatively, the same reasoning applies to the two proposed models discussed in the next section.

## 3. Proposed models

### 3.1 Mean-field interaction-based model

The MF model is the first model proposed in this study, which accounts for the interactions between each $X(r,b\pi(\mathrm{d}r))$ using its weighted average in the jump part. Specifically, given $w\in[0,1]$, the MF model is formulated for $t\in\mathbb{R}$ as

$$Z_t = \int_{r=0}^{r=+\infty} \left\{ \frac{b\pi(\mathrm{d}r)}{r} + \int_{s=-\infty}^{s=t} \int_{z=0}^{z=+\infty} e^{-r(t-s)} z N_{\mathrm{MF}}(\mathrm{d}s,\mathrm{d}z,\mathrm{d}r) \right\}. \qquad (9)$$

The jump intensity of the Poisson random measure $N_{\mathrm{MF}}$ is given as

$$\left\{ wrX_{s-}(r,b\pi(\mathrm{d}r)) + (1-w)\mathbb{E}\left[\int_{y=0}^{y=+\infty} yX_{s-}(y,b\pi(\mathrm{d}y))\right]\pi(\mathrm{d}r)\right\}\nu(\mathrm{d}z)\mathrm{d}s. \qquad (10)$$

The coefficient $w$ is the strength of the terms with and without expectation, and model (9) reduces to the previous one (4) when $w=1$. The mean-field term $\mathbb{E}\left[\int_{y=0}^{y=+\infty} yX_{s-}(y,b\pi(\mathrm{d}y))\right]$ is the expectation of the weighted average of $X$. Weighting by $y$ corresponds to the scaling of $X$ by the reversion speed, analogous to the self-exciting term $wrX_{s-}(r,b\pi(\mathrm{d}r))$. This scaling leads to the decomposition of each moment into quantities on different timescales. Here, the expectation $\mathbb{E}[\cdot]$ can be understood as the operation $\lim_{T\to+\infty}\frac{1}{T}\int_0^T (\cdot)\mathrm{d}t$ by assuming ergodicity as in conventional superposed processes. Thus, we understand that the MF model assumes that the nature of the mean field results from the average of the past states of the system.

The model (9) is seemingly more complex than model (4); however, this issue is resolved by the ansatz

$$\mathbb{E}\left[\int_{y=0}^{y=+\infty} yX_t(y,b\pi(\mathrm{d}y))\right] = Cb, \quad t\in\mathbb{R} \qquad (11)$$

with a constant $C$. Then, the MF model (9) becomes the previous one (4) with some shift, where the nominal process $X(r,b)$ now satisfies the following equation instead of (2):

$$X_t(r,b) = \frac{b}{r} + \int_{s=-\infty}^{s=t} \int_{z=0}^{z=+\infty} \int_{u=0}^{u=rwX_{s-}(r,b)+(1-w)Cb} e^{-r(t-s)} z N(\mathrm{d}s,\mathrm{d}z,\mathrm{d}u), \quad t\in\mathbb{R}. \qquad (12)$$

The well-posedness of (12) follows from Proposition 5.6 of Fu and Li [58], and $C$ in (11) is found below.

*Proposition 1*

$$C = \frac{1}{1-M_1} > 0. \qquad (13)$$

By **Proposition 1**, the average, variance, and autocorrelation of $X$ in (12) are given by

$$\mathbb{E}[X_t] = \frac{b}{1-M_1}\frac{1}{r}, \quad \mathbb{V}[X_t] = \frac{M_2 b}{2(1-M_1)(1-wM_1)}\frac{1}{r}, \quad \rho(\tau) = e^{-r(1-wM_1)\tau}, \qquad (14)$$

respectively. Higher-order statistics, such as the skewness $\mathbb{S}[X_t]$, is also analytically obtained:

$$\mathbb{S}[X_t] = (\mathbb{V}[X_t])^{-3/2} \frac{b}{(1-M_1)(1-wM_1)}\left(\frac{1}{3}M_3 + \frac{1}{2}wM_2^2\right)\frac{1}{r}. \qquad (15)$$

We then obtain the average, variance, autocorrelation, and skewness of the superposed process $Z$ of the MF model (**Section A3**):

$$\mathbb{E}[Z_t] = \frac{b}{1-M_1}R, \quad \mathbb{V}[Z_t] = \frac{M_2 b}{2(1-M_1)(1-wM_1)}R, \quad \rho(\tau) = \frac{1}{R}\int_0^{+\infty}\frac{1}{r}e^{-r(1-wM_1)\tau}\pi(\mathrm{d}r), \qquad (16)$$

and

$$\mathbb{S}[Z_t] = (\mathbb{V}[Z_t])^{-3/2} \frac{b}{(1-M_1)(1-wM_1)}\left(\frac{1}{3}M_3 + \frac{1}{2}wM_2^2\right)R. \qquad (17)$$

The total number of jumps $\mathbb{J}[Z_t]$ of $Z$ within a time interval is again given by (8).

The key similarities and differences between the statistics of (4) and the MF models are as follows: first, they share the same average and total number of jumps per unit time. Hence, the mean-field interaction in the assumed form does not affect the mean behavior of the superposed process. Conversely, their variances and autocovariances differ; hence, the mean-field interaction affects the fluctuation of the superposed process. According to (16), the mean-field interaction ($w<1$) decreases the variance. In the jump rate of the MF model, the random component is $wrX_{s-}(r, b\pi(\mathrm{d}r))$, which becomes less dominated as $w$ decreases, leading to a reduced variance. In addition, the existence of the mean-field interaction increases the exponent $r(1-wM_1)$ in the integrand of the autocorrelation, suggesting that the self-exciting nature of the superposed process $Z$ weakens as the mean-field effect dominates.

### 3.2 Aggregation-based model

We also proposed the AG model, which accounts for interactions among each $X(r, b\pi(\mathrm{d}r))$ through their weighted average in the jump part. More specifically, given $w \in [0,1]$, the superposed process $Z$ is formulated for $t \in \mathbb{R}$ as

$$Z_t = \int_{r=0}^{r=+\infty}\left\{\frac{b\pi(\mathrm{d}r)}{r} + \int_{s=-\infty}^{s=t}\int_{z=0}^{z=+\infty}e^{-r(t-s)}zN_{\mathrm{AG}}(\mathrm{d}s,\mathrm{d}z,\mathrm{d}r)\right\}. \qquad (18)$$

The jump intensity of the Poisson random measure $N_{\mathrm{AG}}$ is given by

$$\left\{ wrX_{s-}\left(r,b\pi(\mathrm{d}r)\right)+(1-w)\int_{y=0}^{y=+\infty} yX_{s-}\left(y,b\pi(\mathrm{d}y)\right)\pi(\mathrm{d}r)\right\}\nu(\mathrm{d}z)\mathrm{d}s. \tag{19}$$

The absence of expectation $\mathbb{E}$ is the difference between the MF and AG models.

The AG model cannot be simplified like the MF model in **Section 3.1** because the integral in (19) is not a constant but a stochastic process. Therefore, we encounter a more difficult case than in the AG model. Nevertheless, the finite-dimensional version (**Section A1**) gives a hint to compute the statistics of $Z$ through the moment-generating function $\mathbb{M}_Z(\theta)=\mathbb{E}\left[e^{-\theta Z_t}\right]$ defined for $\theta \geq 0$. In contrast to the MF model, the absence of expectation $\mathbb{E}$ in the jump rate of the AG model implies that the interaction assumed in this model is only based on the current system and not its time average. Hence, the AG model is more transient. Moreover, as demonstrated later, this difference between both proposed models becomes more visible in their variances.

At a stationary state, for any $\theta \geq 0$, the moment-generating function of $Z$ should be

$$\mathbb{M}_Z(\theta)=\mathbb{E}\left[e^{-\theta Z_t}\right]=\exp(-A(\theta)) \tag{20}$$

with

$$A(\theta)=b\int_{t=0}^{t=+\infty}\int_{r=0}^{r=+\infty}B_t(r,\theta)\pi(\mathrm{d}r)\mathrm{d}t, \tag{21}$$

where $B_t(r,\theta)$ ($t>0$ and a.s. with respect to $\pi$) formally solves the generalized Riccati equation, which is an infinite-dimensional version of that in **Section A1**:

$$\frac{\partial}{\partial t}B_t(r,\theta)=\underbrace{r}_{\text{Time scale}}\times\left\{\begin{array}{l}\underbrace{-B_t(r,\theta)}_{\text{Reversion}}+\underbrace{w\int_0^{+\infty}\left(1-e^{-B_t(r,\theta)z}\right)\nu(\mathrm{d}z)}_{\text{Self-excitation}}\\+\underbrace{(1-w)\int_{y=0}^{y=+\infty}\int_{z=0}^{z=+\infty}\left(1-e^{-B_t(y,\theta)z}\right)\nu(\mathrm{d}z)\pi(\mathrm{d}y)}_{\text{Aggregation}}\end{array}\right\} \tag{22}$$

subject to the initial conditions of $B_0(\cdot,\theta)=\theta$. Each term on the right side of (22) has a clear meaning, as indicated above. Furthermore, the timescale, which is the reciprocal of $r$, appears as a proportional coefficient, theoretically suggesting that each moment of the process is a superposition (integration) of quantities with distributed timescales.

The moment-generating function gives the moment ($\mathbb{E}\left[Z_t^k\right]=(-1)^{k-1}\left.\frac{\mathrm{d}^k\mathbb{M}_Z(\theta)}{\mathrm{d}\theta^k}\right|_{\theta=0}$) when they exist, and the main task is analyzing the initial-value problem of the generalized Riccati equation. We postpone this issue after analyzing linearized generalized Riccati equations (Lyapunov equations presented later) because the former equations are easier to handle.

*Remark 3* The generalized Riccati equation is well-posed for all $w\in[0,1]$ (see **Proposition 3**); therefore, it covers the previous model presented in **Section 2** as a special case of $w=1$.

*Remark 4* The generalized Riccati equation is similar to the partial integro-differential equations of structured neural networks and fields (Section 2.5 of Bressliff [60]). Particularly, the case $w=1$, where jumps in $Z$ are only due to the mean-field effect, would be a special case of their equations. This connection between the integro-differential equations is due to neural activities being triggered by certain neural activities.

In application, we may not need the generalized Riccati equation (22) but its sensitivity $E_t(k,r) = \left. \dfrac{d^k B_t(r,\theta)}{d\theta^k} \right|_{\theta=0}$ ($k=1,2,3,\ldots$) due to

$$\mathbb{E}\left[Z_t^k\right] = (-1)^k \left.\frac{d^k \mathbb{M}_Z(\theta)}{d\theta^k}\right|_{\theta=0} = (-1)^k \left.\frac{d^k \exp(-A(\theta))}{d\theta^k}\right|_{\theta=0}. \tag{23}$$

Particularly, based on (21), we have ($A(0)=0$)

$$\mathbb{E}[Z_t] = -\left.\frac{d\exp(-A(\theta))}{d\theta}\right|_{\theta=0} = b\int_{t=0}^{t=+\infty}\int_{r=0}^{r=+\infty} E_t(k,r)\pi(dr)dt \tag{24}$$

and

$$\mathbb{V}\left[Z_t^2\right] = \mathbb{E}\left[Z_t^2\right] - \left(\mathbb{E}[Z_t]\right)^2 = b\int_{t=0}^{t=+\infty}\int_{r=0}^{r=+\infty} E_t(2,r)\pi(dr)dt. \tag{25}$$

The governing equations of $E_t(k,r)$ ($k=1,2,3,\ldots$), the Lyapunov equations, are obtained as follows:

$$\frac{\partial}{\partial t}E_t(1,r) = r\left\{-(1-wM_1)E_t(1,r) + (1-w)M_1\int_{y=0}^{y=+\infty} E_t(1,y)\pi(dy)\right\}, \tag{26}$$

subject to the initial condition $E_0(1,\cdot)=1$, and

$$\frac{\partial}{\partial t}E_t(2,r) = r\left\{-(1-wM_1)E_t(2,r) + (1-w)M_1\int_{y=0}^{y=+\infty} E_t(2,y)\pi(dy)\right\} \\ + rM_2\left\{w\left(E_t(1,r)\right)^2 + (1-w)\int_{y=0}^{y=+\infty}\left(E_t(1,y)\right)^2 \pi(dy)\right\}, \tag{27}$$

subject to the initial condition $E_0(2,\cdot)=0$.

A few more concepts are required to work through these equations. We set a constant $T>0$, the Lebesgue space $\mathbb{L}^1$ of function $\phi:(0,+\infty)\to\mathbb{R}$ equipped with the norm $\|\phi\|_1 = \int_0^{+\infty}|\phi(r)|\pi(dr)$, and the Banach space $\mathbb{L}_T^1$ of function $\Phi:[0,T]\to\mathbb{L}^1$ equipped with the norm $\|\Phi\|_{1,T} = \sup_{0\le t\le T}\|\Phi_t\|_1$. For each $\lambda > 0$, set another norm $\|\cdot\|_{1,T,\lambda}$ equivalent to $\|\cdot\|_{1,T}$ as $\|\Phi\|_{1,T,\lambda} = \sup_{0\le t\le T}\left(e^{-\lambda t}\|\Phi_t\|_1\right)$.

**Proposition 2** reveals that the Lyapunov equations are well-posed.

*Proposition 2*

*The integro-differential equation (26) subject to the initial condition $E_0(1,\cdot)=1$ admits a unique solution in $\mathbb{L}_T^1$ for any $T>0$ and is continuous with respect to time $t\ge 0$. The same applies to the integro-differential equation (27) subject to the initial condition $E_0(2,\cdot)=0$.*

We also have the following result regarding the generalized Riccati equation. The proof of **Proposition 3** is analogous to that of **Proposition 2** but with a technical modification to deal with the nonlinearity.

*Proposition 3*

*The generalized Riccati equation (22) subject to an initial condition $B_0(\cdot, \theta) = \theta$ ($\theta \geq 0$) admits a unique nonnegative solution in $\mathbb{L}_T^1$ for any $T > 0$ and is continuous with respect to time $t \geq 0$.*

More explicit results can be obtained by assuming that the last integrals in (24) and (25) exist. More specifically, we have (**Section A4**)

$$\mathbb{E}[Z_t] = \frac{b}{1-M_1} R, \quad \mathbb{V}[Z_t] = \frac{bM_2}{(1-M_1)(1-wM_1)} \left\{ \frac{1}{2} R + (1-w) M_1 \int_{t=0}^{t=+\infty} \left( \int_{r=0}^{r=+\infty} E_t(1,r) \pi(\mathrm{d}r) \right)^2 \mathrm{d}t \right\}. \quad (28)$$

Hence, in practice only $E_t(1,r)$ and its integral are necessary for computing the average and variance. The left side of equation (28) demonstrates that the previous, MF and AG models, share the same average, and we can obtain the same number of jumps per unit time: $\mathbb{J}[Z_t]$. Therefore, the three models have the same low-order statistics. Furthermore, the variance differs, as indicated on the right side of equation (28). According to (5), (14), and (28) the variance in the AG model is larger than that in the MF model. Hence, the MF model has the smallest fluctuations among the three models. We found no analytical results for the order of the previous and AG models. Therefore, we computationally investigated this issue. We are unaware of the representation formula of autocorrelation, which will also be studied computationally.

## 4. Application
### 4.1 Target species

The superposed processes are applied to the biological count data of migrating fish—*P. altivelis*—in five rivers in Japan. Each river flows into the Pacific Ocean, serving as a migratory route for *P. altivelis*. We focused on their spring upstream migration from the Pacific Ocean to a river, which occurs between March and July yearly. The migration duration and size differ among rivers and years. Investigating the upstream migration of *P. altivelis* from spring to summer is important because it controls their population dynamics in rivers from summer to autumn. This population dynamics affect trophic cascades in and around rivers [61,62], of which this is a vital commercial species in inland fisheries in Japan [e.g., 63,64]. A biological count should be an integer; however, we regarded it as a continuous variable by assuming that the count is significantly larger than one, on average.

### 4.2 Target rivers

The target rivers in this study are the Yahagi, Nagara, Tama, Ara, and Tone Rivers, into which *P. altivelis*

migrates yearly (**Figure 1**). **Table 1** presents the daily count data for individual *P. altivelis*. The stochastic analyses of the data of the Tama, Tone, and Ara Rivers have not been addressed to the best of the author's knowledge. These data would not include all migrating individual fish but only a part of them for technical reasons; nonetheless, they would represent *P. altivelis* migration.

Upstream migration of *P. altivelis* has recently been investigated in each river as follows: Yahagi River (fish counts by Yamamoto and Nagatomo [65] and Yamamoto et al. [66]), Nagara River (mass balance model of the fish population dynamics by Mouri et al. [67]), Tama River (fish count methodologies by Takase et al. [68]), Ara River (swimming ability by Shiina et al. [69]), and Tone River (influences of a weir on the upstream migration by [70]). Some Yahagi River data have been studied using only the supOU process [51] and a system of nonlinear hybrid SDEs [71]. However, these data have not been analyzed using a unified method like the proposed ones, which motivated us to apply superposed processes. The stationarity assumption of the present superposed processes should be removed in the future because fish migration is part of the seasonal dynamics of *P. altivelis*. Nevertheless, the application described in this section provides the primary result for a more complex and advanced analysis.

All the data used in this study are the daily counts of upstream-migrating fish at fixed observation points. The data were curated in a unified manner. We curated the data such that the first and last zeros in each time series were omitted for each observation station and year. $I$ denotes the total number of data points. The first to last days are numbered Day 1 to Day $I$. These procedures were applied to all the datasets. **Table S1** reports the $I$ and Day 1 for each dataset. **Figure 2** presents examples of the time series data for the Nagara River, suggesting that each time series has several spikes. Other data for each river are numerically obtained from each open source (**Table 1**).

**Table S2** summarizes the empirical average (Ave), standard deviation, coefficient of variation (CV), total number of jumps per unit time (Jmp: total number of interior and boundary strict local maximum points in time series data, assuming that no strict maxima exists for unavailable), and skewness (Skw) of each data point. According to **Table S2**, the average and variance depend on the year and river and vary by several orders of magnitude in one river. In contrast, in all the time series data, the empirical CV has the order $O(10^0)$, with the mean of all cases being 1.82. The total number of jumps in a unit time Jmp has the order of $O(10^{-1})$, with the mean of all cases being 0.283, revealing that one jump event occurs every three to four days on average. All the time series data are positively skewed in the order of $O(10^0)$, with the mean of all cases being 2.89. Hence, their distributions are positively skewed, suggesting that a suitable jump can capture them. The common magnitudes of CV, Jmp, and skewness and the common sign of skewness among the data imply that these statistics effectively characterize the biological count data of *P. altivelis*. The CV values imply that the average and variance of the daily fish counts are related.

**Table 1.** Sources of the *P. altivelis* count data.

| River | Years | Sources | Remarks |
|---|---|---|---|
| Yahagi River | 1998-2023 | 1998 to 2009: Yamamoto and Nagatomo [65]<br>2010 to 2020: Yamamoto et al. [66]<br>2021-2022: Natural Ayu Ecological Survey Executive Committee [72]<br>2023: Yahagi River Fisheries Cooperative [73] | ✓ Data was mainly collected at the left bank fishway of Meiji-Yousui Irrigation Head Works 35 (km) upstream from the river mouse.<br>✓ The data in 2016, 2020, 2021 were collected at the right bank fishway.<br>✓ Daily count data in 2024 are not available.<br>✓ Used video counting system. |
| Nagara River | 2023-2024 | Nagaragawa Estuary Barrage Operating & Maintenance Office [74,75] | ✓ Data were collected at the Nagaragawa Estuary Barrage 5.4 (km) upstream from the river mouse.<br>✓ Used video counting system. |
| Tama River | 2011-2024 | Tokyo Metropolitan Islands Area Research and development Center for Agriculture, Forestry, and Fisheries [76]. | ✓ Data were collected at 11 (km) upstream from the river mouse.<br>✓ Used fixed fishing net. |
| Ara River | 2012-2024 | Tone Water Direction General Management Office [77]. | ✓ Data were collected at the Akigase Intake Weir 35 (km) upstream from the river mouse.<br>✓ Used fixed fixing net. |
| Tone River | 2012-2024 | Tone Water Direction General Management Office [77]. | ✓ Data were collected at the Tone Great Weir 154 (km) upstream from the river mouse.<br>✓ Used fixed fishing net. |

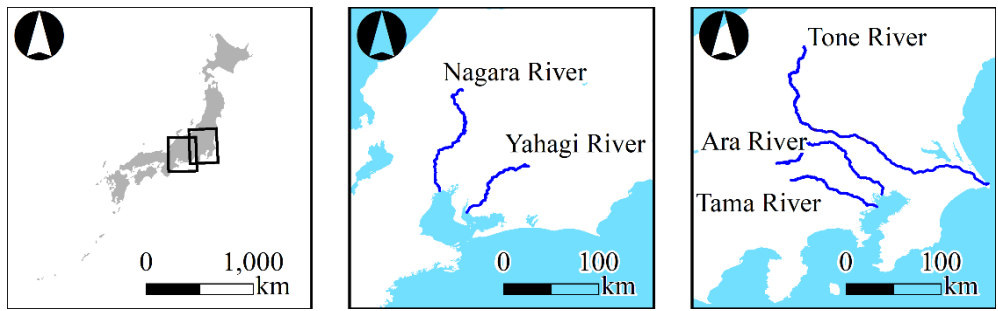

**Figure 1.** Map of the studied rivers.

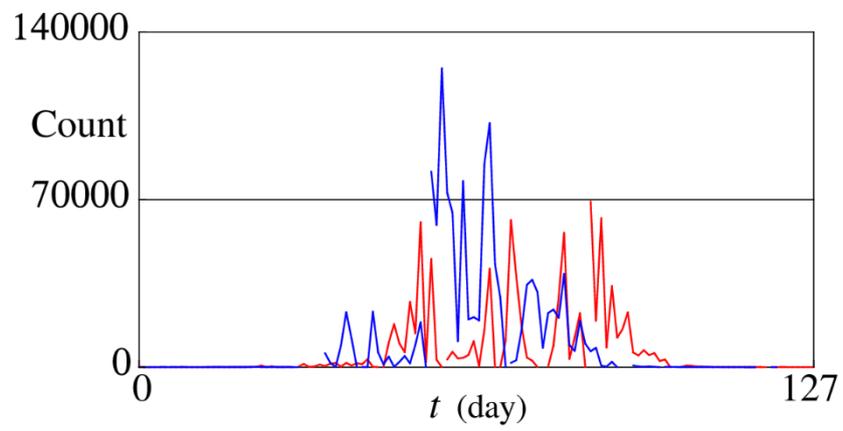

**Figure 2.** Time series data of the fish count at the Nagara River in 2023 (red) and 2024 (blue).

We examined the relationships between these statistics and concluded a high correlation between CV and skewness (both non-dimensional), which is Skw $= \kappa$CV with $\kappa = 1.5875$ ($R^2$=0.95) (**Figure 3**). This implies that the time series with a higher fluctuation is more positively skewed. In the MF model, this suggests

$$\frac{\mathbb{S}[Z_t]}{(\mathbb{V}[Z_t])^{3/2}} = \kappa \frac{(\mathbb{V}[Z_t])^{1/2}}{\mathbb{E}[Z_t]} \tag{29}$$

or equivalently

$$\frac{1}{3}M_3 + \frac{1}{2}wM_2^2 = \kappa \frac{M_2^2}{4(1-wM_1)}. \tag{30}$$

Furthermore, we assume that the Lévy measure $v$ has the exponential density $\mu \times \lambda e^{-\lambda z}$ ($z > 0$), with the jump frequency $\mu > 0$ and jump size $\lambda^{-1} > 0$. This Lévy measure $v$ is the simplest model of background positive-jump processes. In this case, based on (30), we have the quadratic equation of $M_1 = \mu / \lambda$:

$$w^2 M_1^2 + \frac{\kappa}{2} M_1 - 1 = 0. \tag{31}$$

If $w = 0$, then $M_1 = \frac{2}{\kappa} = 1.216 > 1$, which violates the assumption $M_1 \in (0,1)$. If $w > 0$, then $M_1 = \frac{1}{4w^2}\left(-\kappa + \sqrt{\kappa^2 + 16w^2}\right)$ (positive solution). This $M_1$ belong to $(0,1)$ if $w > \sqrt{1-\frac{\kappa}{2}} = 0.422$. This implies that a model with a small $w$ may not be feasible. This rough yet theoretical estimate is supported by the realizability condition of the model parameters discussed in the next subsection, which is derived using a different approach than that employed here.

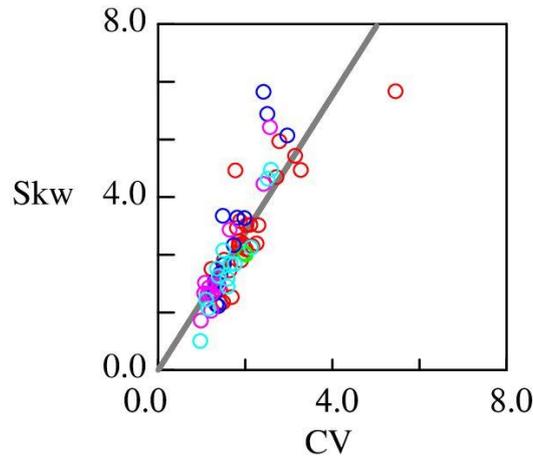

**Figure 3.** The relationship between the skewness (Skw) and coefficient of variation (CV). Colors are used to indicate the rivers: Yahagi River (red), Nagara River (green), Tama River (blue), Ara River (magenta), and Tone River (light blue). The grey line represents the linear fit Skw $= \kappa$CV.

## 4.3 Model fitting

We fit each time series data point as a superposed process denoted as $Z$. We fit the MF model using a least-squares estimation of autocorrelation, along with a moment-matching method [34]. We used the Lévy measure with an exponential density $\mu \times \lambda e^{-\lambda z}$. We assume that the probability measure $\pi$ is the Gamma distribution with a density proportional to $r^{\alpha-1} e^{-r/\beta}$ ($r > 0$), with shape parameter $\alpha > 1$ and scale parameter $\beta > 0$, as in a previous study (e.g., Yoshioka et al. [51]). Subsequently, we set $\tilde{\beta} = \beta(1 - wM_1)$. First, we assume that $w \in [0,1]$ is given and subsequently examine cases where $w$ is also fitted. The autocorrelation with lag $\tau \geq 0$ is then theoretically given by $\rho(\tau) = (1 + \tilde{\beta}\tau)^{-(\alpha-1)}$, and the corresponding superposed process has a truly long memory ($\int_0^{+\infty} \rho(\tau) d\tau = +\infty$) if and only if $1 < \alpha \leq 2$. This autocorrelation also covers the other extreme case with exponential decay because $\rho(\tau) = (1 + \tilde{\beta}\tau)^{-(\alpha-1)} \to \exp(-\gamma\tau)$, under the limit $\alpha \to +\infty$ (or $\beta \to 0$), where we assume the constant $\gamma > 0$ exists, such that $\gamma = (\alpha - 1)\tilde{\beta}$. Moreover, with this autocorrelation, the Hurst exponent $H$ is given by $H = \frac{3}{2} - \frac{1}{2}\alpha$, if it is positive and between 0 and 1 and hence if $\alpha \in (1,3)$ [e.g., 78,79] due to $\rho(\tau) \sim \tau^{-(\alpha-1)}$ for large $\tau$. Truly long memory cases correspond to $H \in \left(0, \frac{1}{2}\right)$. Consequently, the Hurst exponent of the fish count data can be estimated by fitting the MF model against them. This is an advantageous point of the MF model in practice.

In the model fitting, the parameters $\alpha, \tilde{\beta}$ are first identified via a classical nonlinear least-squares fitting between the empirical and theoretical autocorrelations against the time lags $1, 2, 3, ..., 14$ (days) because each jump seems to continue empirically for a maximum of 14 days. Moment fitting was then applied between the empirical and theoretical Ave, Var, and Jmp as follows: $\mathbb{E}[Z_t] = \text{Ave}$, $\mathbb{V}[Z_t] = \text{Var}$, and $\mathbb{J}[Z_t] = \text{Jmp}$. These equations can be rewritten to obtain $\lambda$, $b$, and $\mu$ in this order, where $M_0 = \mu$, $M_1 = \mu\lambda^{-1}$, $M_2 = 2\mu\lambda^{-2}$, and $\beta = \tilde{\beta}(1 - wM_1)^{-1}$:

$$\lambda = \sqrt{\frac{1}{\tilde{\beta}(\alpha - 1)} \frac{\text{Jmp}}{\text{Var}}} > 0, \tag{32}$$

$$b = \text{Ave}\tilde{\beta}(\alpha - 1) - (1 - w)\frac{\text{Jmp}}{\lambda} \in \mathbb{R}, \tag{33}$$

$$\mu = \frac{\text{Jmp}}{\text{Jmp} + b\lambda} \lambda \in (0,1). \tag{34}$$

In principle, no modeling errors existed in the average or variance of this fitting. The $\lambda$ estimate is positive, and that of $\mu$ is between 0 and 1, as desired. In contrast, the $b$ estimate becomes positive only if

$$1 \geq w > \max\left\{0, 1 - \sqrt{\frac{\tilde{\beta}(\alpha-1)}{\text{Jmp} \cdot \text{CV}^2}}\right\}. \tag{35}$$

The inequality in (35) serves as the realizability constraint of $w$ to guarantee the positivity of the identified value of $b$. Therefore, we cannot assume $w$ to be close to 0 in some cases. Moreover, the theoretical estimate of $M_1$ in the previous subsection also suggests that a small $w$ may give unrealistic results.

### 4.4 Comparison among models
#### 4.4.1 MF model

Initially, we assume $w=1$, in which condition (35) is always satisfied. **Table S3** summarizes the values of the identified parameter at each data point and river. **Figure 4** illustrates the empirical and fitted ACFs for the selected cases with long or exponential memory. The identified values of the shape parameter $\alpha$ of the probability measure $\pi$ suggest that the memory lengths of fish count data vary with the years, even in one river. For instance, in the Yahagi River, the order of $\alpha$ is from $O(10^0)$, including those with truly long memory, to $O(10^4)$, including those with exponential memory. Qualitatively the same results apply to the Tama and Ara Rivers. In the Nagara River, the data revealed a truly long memory. The Tone River have $\alpha$ larger than $O(10^1)$; hence, fish migration in this river is suggested not to persist and is impulsive. Applying a model without a long memory, such as the classical Ornstein–Uhlenbeck process or its self-exciting version, would suffice for studying the upstream migration of *P. altivelis* in this river. Concerning the Hurst exponent $H$, 30% data of the Yahagi River, the two data of the Nagara River, and more than the half data of the Tama River show $H \in (0,1)$. The estimated Hurst exponents take values larger than $1/2$ for the Nagara River, and the Yahagi and Tama have both cases $H > 1/2$ and $H < 1/2$; we did not find explicit criteria to separate the two cases.

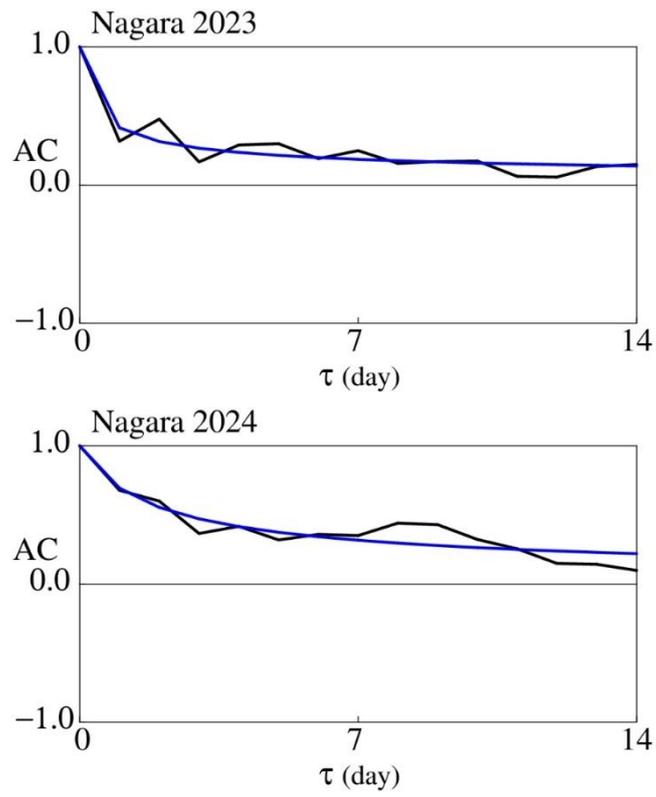

**Figure 4.** Autocorrelation (AC) for the Nagara River data in 2023 and 2024: empirical (black) and theoretical results (blue).

We investigate cases where the parameter $w$ is also fitted. The fitting procedure is the same as in the previous case with $w$ being fixed except that it fitted by solving the least-squares problem to minimize the relative error between the empirical and theoretical skewness. The fitted value of $\alpha$ and the memory do not change owing to this fitting procedure modification. **Table S4** summarizes the fitted parameter values, and **Table S5** shows a comparison of the theoretical skewness values and their relative errors between the cases with and without fixing $w$. Additionally, no modeling errors exist in the average and variance. Furthermore, adding $w$ to the set of parameters to be fitted improves the reproductivity of skewness, e.g., more than 80% and 40 % decrease of the Skw prediction error for the Nagara River in 2023 and 2024, respectively; however, the optimized $b$ values become negative in several cases, as indicated by magenta in **Table S5**. This implies that fitting $w$ using the skewness, a high-order statistic, would not necessarily guarantee the realizability of the model.

Most of the data suggest $w = 0$ for the Ara and Tone Rivers, except for the few cases, suggesting that fish migration is not due to self-exciting jumps but is caused by external factors such as temperature changes in the sea, lake, or river and river flows [67]. In contrast, Sato and Seguchi [80] reported no significant correlation between the number of fish in a fishway, water temperature, and salinity, suggesting a contrasting result and the complexity of this biological phenomenon. This finding, combined with the short memory of the data of the Tone River, suggests that the classical Ornstein–Uhlenbeck process could be used to study fish migration in this river. Intermediate $w$ values were observed in some rivers, implying the coexistence of self- and non-self-exciting jumps.

The realizability condition in (33) reveals that a small $w$ may yield an unrealistic model with a nonpositive $b$. Such cases occurred in all rivers except the Ara River. This implies that using higher-order statistics, such as skewness, in identifying parameters is not always effective. However, **Table S5** suggests that it is effective without degrading the reproducibility of the lower-order statistics, such as the average and variance. No absolute criterion exists for matching the moment; thus, we recommend examining methods with and without using higher-order statistics and comparing their performances, as in this study.

We further explored the link between memory and jumps through non-dimensionalization. An elementary calculation revealed that the following quantities have units of time and individuals (ind):

$$\bar{t} = \frac{R}{1-M_1} = \frac{1}{1-M_1}\frac{1}{\beta(\alpha-1)} \text{ and } \bar{X} = \frac{bR}{1-M_1} = \frac{b}{1-M_1}\frac{1}{\beta(\alpha-1)}. \qquad (36)$$

We then used the following non-dimensionalization that corresponds to the model in which $b$ is reduced to 1 and $\mathbb{E}[Z_t]$ to 1: $t \to t\bar{t}$, $X_t \to X_t\bar{X}$, $\lambda \to \lambda/\bar{X}$, $\mu \to \mu/\bar{X}$, $z \to z\bar{X}$, $r \to r/\bar{t}$, $\beta \to \beta/\bar{t}$, $\alpha \to \alpha$, and $w \to w$. In the sequel, we deal with non-dimensionalized models.

We explored the link between $\beta$, as the (reciprocal of) timescale, and $M_1 = \mu/\lambda$, as the contribution of jumps. $\beta$ varies across scales (i.e., about eight orders); thus, we explore the relationship between $M_1$ and $\text{Log}\beta$ with "Log" as the ordinary logarithm (**Figure 5**). We examine the case of $w$ being fixed at one and that of $w$ being fitted, and the data with negative $b$ were excluded from the latter. **Figure**

**5(a)** reveals that the plots for $w=1$ are scattered over the panel; however, they are not uniformly distributed but are separated into left and right clusters with small and large $\beta$, respectively. Each river, except for the Nagara River, which has only two data points, contributes to both clusters. Even in one river, the timescale and jump contribution to fish migration dynamics differ with the years. A similar finding applies to **Figure 5(b)** with $w$ being fitted but with a wider range with respect to $M_1$. $\beta$ and $M_1$ are essential in characterizing the superposed process; therefore, the MF model suggests that fish migration has two phases in terms of timescale and jump contribution to dynamics.

The MF model is conceptual and does not account for the detailed mechanisms of migration triggers and factors. Nonetheless, the model-fitting results suggest that studying fish count data based on a stochastic process model would aid the comprehension of this complex biological phenomenon. Particularly, the results of this study suggest that an exponential memory or a long memory model is insufficient; however, a unified method that can cover both is necessary to analyze fish migration.

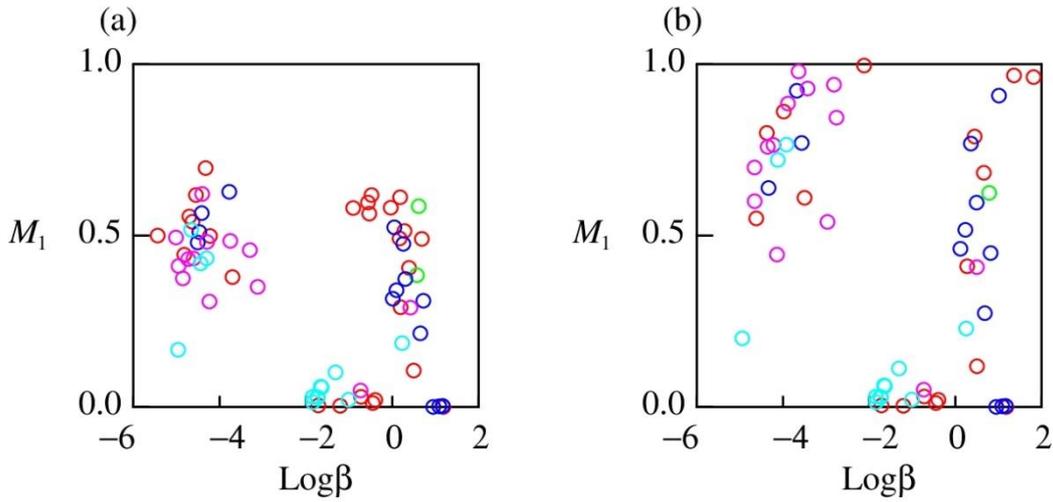

**Figure 5.** The relationship between the normalized $\beta$ and $M_1$: **(a)** $w=1$ and **(b)** $w$ being fitted. Colors are used to indicate the rivers: Yahagi River (red), Nagara River (green), Tama River (blue), Ara River (magenta), and Tone River (light blue).

**4.4.2 AG model**

We compare the performance of both proposed models numerically because the AG model does not have closed-form statistics other than the average. We continue using the non-dimensionalized system; hence, the parameter $b$ is fixed to 1 and the average of $Z$ to 1. We investigate truly long ($\alpha \in (1, 2)$), critical ($\alpha = 2$), and moderately long memory cases ($\alpha > 2$).

We computed the variance using the generalized Riccati equation discretized under a naïve finite difference method [81] and autocorrelation using the Monte Carlo method [82]. The finite difference and Monte Carlo methods compute the Riccati-type equations and supOU and related non-Markovian SDEs, respectively. Both methods use the finite-dimensional version of the superposed processes and their associated generalized Riccati equations presented **in Section A1**. The Lyapunov equations and the Monte Carlo method are used to compute the variance and the memory, respectively. The time increments for temporal discretization are 0.001 (days) and 0.0005 (days) for the finite difference and Monte Carlo methods, respectively. The degree of freedom of the finite-dimensional superposed processes are 65,536 and 2,048 for the finite difference and Monte Carlo methods, respectively, considering computational accuracy.

**Figure 6** illustrates the computed sample paths of the AG model for different values of $\alpha$, where we cover a moderately long memory case ($\alpha = 4$), a critical case ($\alpha = 2$), and a truly long memory case ($\alpha = 1.8$). **Figure 6** illustrates the longer memory of the truly long memory case. **Figure 7** presents a comparison of the variances between the MF and AG models. The variance is an increasing function of the weight $w$ for all cases, suggesting that for both models, the self-exciting part $wX$ of jump intensity contributes more to the variance than the mean-field or the aggregation part for a larger $w$. Additionally, the MF model has a higher variance than the AG model, as theoretically suggested in **Section 3.2**. In this case, the variance of the previous model equals 2, which is the largest variance among all the models. **Figure 7** also suggests the influence of memory in the AG model, where longer memory and hence a smaller $\alpha$ result in a variance closer to that of the MF model. This suggests that the prevalence of long memory in both proposed models is somewhat similar; however, it is beyond the scope of this study because a more advanced theoretical stochastic calculus is needed to address this issue.

Finally, **Figure 8** compares the autocorrelations of the AG models for different $\alpha$ and $w$ values, where those for $w < 1$ are computed using the Monte Carlo simulation. Regarding the cases of $\alpha$ covering truly long, critical, and moderately long memory, the tail of autocorrelation becomes heavier as $w$ increases, i.e., as the self-exciting part dominates the aggregation part. This tendency is qualitatively the same as that of the MF model plotted in **Figure 8** panels. The difference between the AG and MF models becomes more visible with smaller $w$, such that the lightening of the tail of autocorrelation becomes more significant in the MF model. Consequently, the comparison between the AG and MF models suggests that the latter is more sensitive to parameter changes than the former. From an engineering perspective, this implies that model uncertainty, such as a misspecification in the probability measure $\pi$ [e.g., 26,83], will affect the MF model considerably. The AG model is more robust, although its identification seems more challenging than that of the MF model, owing to its low tractability. Thus, a unified identification method for superposed processes is required to address this issue in the future.

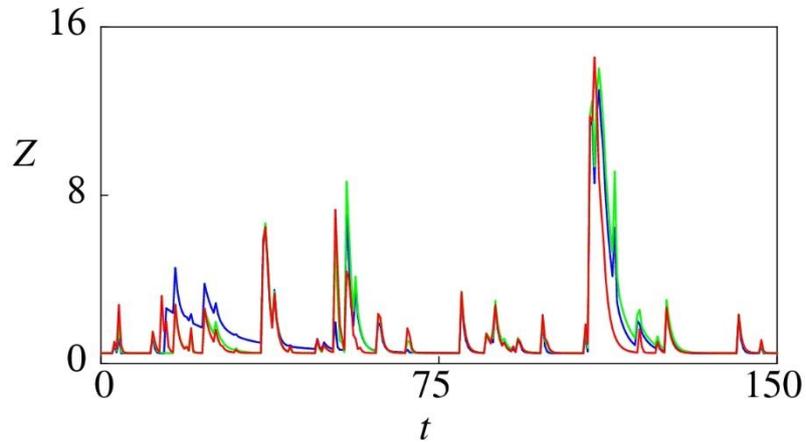

**Figure 6.** Sample paths of the AG model: $\alpha = 4$ (red), $\alpha = 2$ (green), and $\alpha = 1.8$ (black).

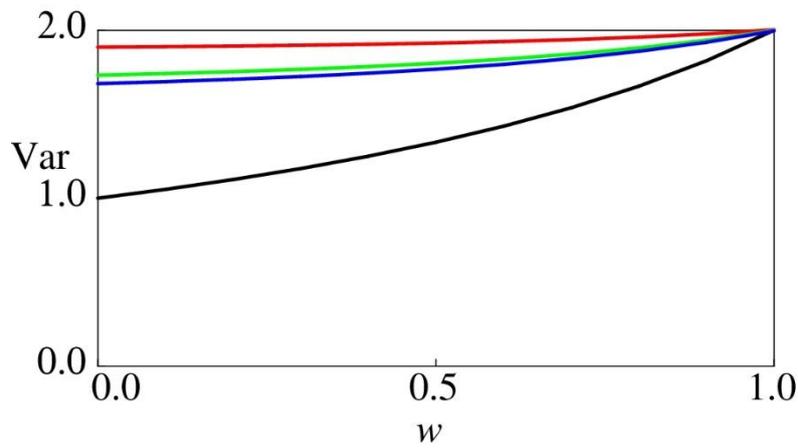

**Figure 7.** Variances (Var) of the MF model (black) and those of AG models ($\alpha = 4$ (red), $\alpha = 2$ (green), and $\alpha = 1.8$ (blue)). The Var of the MF model follows the black curve for all $\alpha > 1$ under the non-dimensionalization.

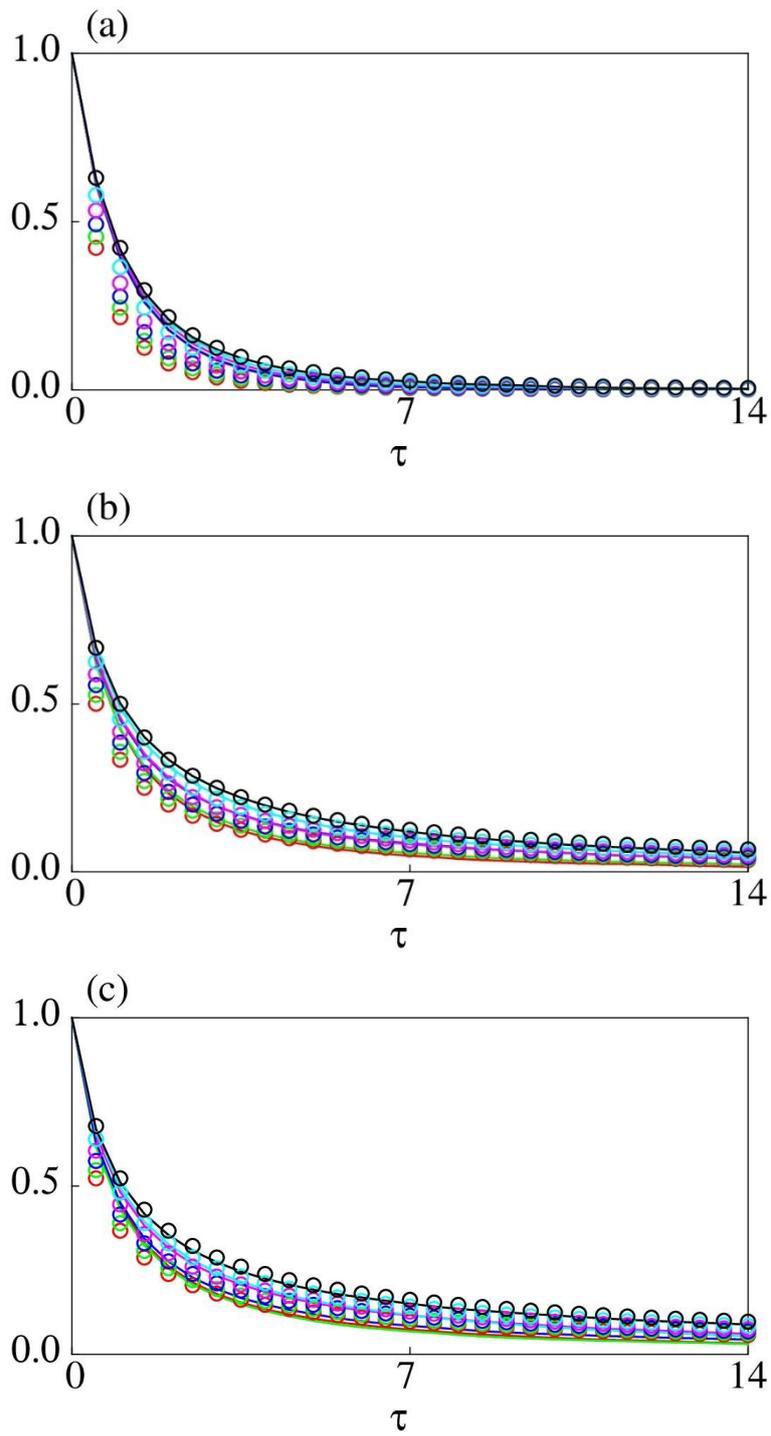

**Figure 8.** Autocorrelations for **(a)** $\alpha = 4$, **(b)** $\alpha = 2$, and **(c)** $\alpha = 1.8$. The lines and circles represent AG and MF models. Colors indicate values as follows: $w = 1$ (black), $w = 0.8$ (light blue), $w = 0.6$ (magenta), $w = 0.4$ (blue), $w = 0.2$ (green), and $w = 0$ (red).

## 5. Conclusion

Long memory processes based on the superposition approach were formulated and analyzed through their mutual interaction via a mean-field or an aggregation term. The statistics of the MF and AG models were discussed. Explicit formulae for the statistics of these models were also derived. An application study to characterize the migrating fish count data of *P. altivelis* in five rivers in Japan suggested that the memory structure differs each year, even in the same river; conversely, they are similar in the order of variation coefficient and the skewness sign. In addition, we numerically compared the proposed models. The results of this study suggest that an exponential memory or a long memory model is insufficient; however, a unified method that can cover both, like our model, is necessary to analyze fish migration.

In theory, a limitation of the proposed models is that only a few interactions are considered. For instance, theoretically, both the mean-field and aggregation effects can be described using a more complex network structure among the processes to be superposed, which in our case may be a graphon [84]. Specifically, a graphon with an unbounded domain differs from classical graphons with bounded domains. This generalization is possible; however, determining the network structure from the data is a potential challenge. A weakness of the proposed models in terms of applications is that they are conceptual, with their social cues considered abstract, whereas a more reasonable model would use some mechanistic interactions among individuals, as in ecological models [85,86]. A conceptual model can be obtained by upscaling an individual-based model, in which the memory difference compared to the latter is essential. Spatio-temporal random fluctuations also become important in ecological dynamics such as fish migration in some applications, which can be addressed by using a compartment model or an advection-diffusion-reaction model [87,88]. How to efficiently formulate long memory processes in these models will be a crucial point.

Another limitation of this study is that we assumed the stationarity of the fish counts for modeling simplicity. It is more reasonable, but more complex, to use the data as nonstationary. Some nonstationary stochastic models [89] can address this issue; however, combining them with the flexible memory structure of superposed processes is still a challenge to address in future studies. The computational cost of simulating superposition processes is also crucial. In our case, the simulating processes for small $\alpha$, such as 1.2, were extremely inefficient. An efficient quadrature for discretizing the superposed processes with respect to the reversion speed will be explored in the future.

# Appendix

## A1. Finite-dimensional model

This study presents a finite-dimensional version of the AG model. Those of the other two models were obtained analogously. The dimensionality $n \in \mathbb{N}$ is chosen; let $\{r_i\}_{1 \leq i \leq n}$ be positive, strictly increasing sequence. Let $\{\pi_i\}_{1 \leq i \leq n}$ be a positive sequence satisfying $\sum_{i=1}^{n} \pi_i = 1$. The probability measure $\pi$ is then approximated by $\sum_{i=1}^{n} \pi_i \delta_{r_i}$, where $\delta_r$ is the 1-D Dirac Delta concentrated at point $r$.

The following system of $n$-dimensional self-exciting SDEs is considered:

$$dY_t(r_i) = (b\pi_i - r_i Y_t(r_i))dt + \int_{z=0}^{z=+\infty} \int_{u=0}^{u=wr_i Y_{t-}(r_i)+(1-w)\sum_{j=1}^{n} r_j Y_{t-}(r_j)} zN_i(dt, dz, du), \ t \in \mathbb{R}, \ 1 \leq i \leq n. \tag{37}$$

Here, $N_i$ are mutually independent copies of the Poisson random measure $N$. At a stationary state, the superposed process $Z^{(n)}$ is defined using the finite-dimensional system as the finite sum:

$$Z_t^{(n)} = \sum_{i=1}^{n} Y_t(r_i), \ t \in \mathbb{R}. \tag{38}$$

We consider the moment-generating function of the process $Z^{(n)}$ at a stationary state:

$$\mathbb{M}_{Z^{(n)}}(\theta) = \mathbb{E}\left[e^{-\theta Z_t^{(n)}}\right], \ \theta \geq 0. \tag{39}$$

We use the Kolmogorov equation associated with the moment-generating function (the right side is the infinitesimal generator (e.g., Theorem 1.22 in Øksendal and Sulem [57]):

$$\frac{\partial \vartheta_t}{\partial t} = \sum_{i=1}^{n} (b\pi_i - r_i x_i)\frac{\partial \vartheta_t}{\partial x_i} + \sum_{i=1}^{n} \left(wr_i x_i + (1-w)\pi_i \sum_{j=1}^{n} r_j x_j\right) \int_{z_i=0}^{z_i=+\infty} \Delta\vartheta_{t,i}(z_i) \nu(dz_i) \tag{40}$$

for $t > 0$ and $(x_1, x_2, ..., x_n) \in (0, +\infty)^n$ subject to the initial condition $\vartheta_0 = \exp\left(-\theta \sum_{i=1}^{n} x_i\right)$. Here, $\vartheta_t = \vartheta_t(x_1, x_2, ..., x_n)$ is the conditional expectation

$$\vartheta_t(x_1, x_2, ..., x_n) = \mathbb{E}\left[e^{-\theta Z_t^{(n)}} \Big| (Y_0(r_1), Y_0(r_2), ..., Y_0(r_n)) = (x_1, x_2, ..., x_n)\right] \tag{41}$$

under the ergodicity ansatz $\mathbb{M}_{Z^{(n)}}(\theta) = \lim_{t \to +\infty} \vartheta_t(x_1, x_2, ..., x_n)$, and we used the notation

$$\Delta\vartheta_{t,i}(z_i) = \vartheta_t(x_1, x_2, ..., x_i + z_i, ..., x_n) - \vartheta_t(x_1, x_2, ..., x_n). \tag{42}$$

Assuming the solution of the form with time-dependent functions $A_t^{(n)}$, $B_{t,i}^{(n)}$ ($i = 1, 2, 3, ..., n$)

$$\vartheta_t(x_1, x_2, ..., x_n) = \exp\left(-A_t^{(n)} - \sum_{i=1}^{n} B_{t,i}^{(n)} x_i\right) \tag{43}$$

and substituting (43) into the Kolmogorov equation (40) yields

$$A_t^{(n)} = b \sum_{i=1}^{n} \pi_i B_{t,i}^{(n)} \tag{44}$$

and

$$\frac{\mathrm{d}}{\mathrm{d}t} B_{t,i}^{(n)} = r_i \left\{ -B_{t,i}^{(n)} + w \int_{z_i=0}^{z_i=+\infty} \left(1 - e^{-B_{t,i}^{(n)} z}\right) \nu(\mathrm{d}z_i) + (1-w) \sum_{j=1}^{n} \int_{z_j=0}^{z_j=+\infty} \left(1 - e^{-B_{t,j}^{(n)} z}\right) \nu(\mathrm{d}z_j) \right\} \tag{45}$$

for $t > 0$ and $i = 1, 2, 3, ..., n$. The initial conditions for these equations are $A_0^{(n)} = 0$ and $B_{t,i}^{(n)} = \theta$ ($i = 1, 2, 3, ..., n$). The generalized Riccati equation (45) is a space discretization of (22); therefore, it admits the same structure. The well-posedness of (45) can be proven by following the strategy used in the **Proof of Proposition 2**, where the probability measure $\pi$ is replaced by $\sum_{i=1}^{n} \pi_i \delta_{r_i}$. After the elementary calculation, we can also obtain the stationary average as follows (they can be derived by directly taking the expectation of (37)):

$$\mathbb{E}\left[Z_t^{(n)}\right] = \frac{b}{1 - M_1} R^{(n)}, \tag{46}$$

with $R^{(n)} = \sum_{i=1}^{n} r_i^{-1} \pi_i$. Owing to $R < +\infty$, we expect to have the convergence $R^{(n)} \to R$ as $n \to +\infty$ under moderate conditions. Higher-order moments can also be obtained using Equations (44)–(46) as described in the main text.

Consequently, the superposed process can be considered the infinite-dimensional limit of the finite-dimensional model presented in this section. The link between the finite- and infinite-dimensional models is suggested by their generalized Riccati equations. Moreover, the appearance of $\pi$ in the jump rate of the superposed process in the main text originates from the weight $\pi_i$ of the finite-dimensional version. This correspondence, though formal, helps to understand the superposed processes.

**A2. Statistics of the previous model**

For the average $\mathbb{E}[Z_t]$, we formally have

$$\begin{aligned}
\mathbb{E}[Z_t] &= bR + \int_{r=0}^{r=+\infty} \mathbb{E}\left[M_1 \int_{-\infty}^{t} e^{-r(t-s)} r X_s(r, b\pi(\mathrm{d}r)) \mathrm{d}s\right] \\
&= bR + M_1 \int_{r=0}^{r=+\infty} \int_{-\infty}^{t} e^{-r(t-s)} r \mathbb{E}\left[X_s(r, b\pi(\mathrm{d}r))\right] \mathrm{d}s \\
&= bR + M_1 \int_{r=0}^{r=+\infty} \int_{-\infty}^{t} e^{-r(t-s)} r \frac{b\pi(\mathrm{d}r)}{r(1 - M_1)} \mathrm{d}s \\
&= bR + \frac{M_1 b}{1 - M_1} \int_{r=0}^{r=+\infty} \int_{-\infty}^{t} e^{-r(t-s)} \pi(\mathrm{d}r) \mathrm{d}s \\
&= bR + \frac{M_1 b}{1 - M_1} R \\
&= \frac{b}{1 - M_1} R
\end{aligned} \tag{47}$$

For the variance $\mathbb{V}[Z_t] = \mathbb{E}\left[\left(Z_t - \mathbb{E}[Z_t^2]\right)^2\right]$, because each $X(r, b\pi(\mathrm{d}r))$ is independent, we have

$$\mathbb{V}[Z_t] = \int_0^{+\infty} \mathbb{V}\left[X_t(r, b\pi(\mathrm{d}r))\right]$$
$$= \int_0^{+\infty} \frac{M_2 b\pi(\mathrm{d}r)}{2r(1-M_1)^2} \qquad . \qquad (48)$$
$$= \frac{M_2 b}{2(1-M_1)^2} R$$

The independence among $X(r, b\pi(\mathrm{d}r))$ also applies to autocovariance as follows:

$$\mathbb{C}[Z_{t+\tau}, Z_t] = \int_0^{+\infty} \mathbb{C}\left[X_{t+h}(r, b\pi(\mathrm{d}r)), X_t(r, b\pi(\mathrm{d}r))\right]$$
$$= \int_0^{+\infty} \frac{M_2 b}{2r(1-M_1)^2} e^{-r(1-M_1)\tau} \pi(\mathrm{d}r) \qquad . \qquad (49)$$
$$= \frac{M_2 b}{2(1-M_1)^2} \int_0^{+\infty} \frac{e^{-r(1-M_1)\tau}}{r} \pi(\mathrm{d}r)$$

Consequently, we obtain

$$\rho(\tau) = \frac{\mathbb{C}[Z_{t+\tau}, Z_t]}{\mathbb{V}[Z_t]} = \frac{1}{R} \int_0^{+\infty} \frac{1}{r} e^{-r(1-M_1)\tau} \pi(\mathrm{d}r) . \qquad (50)$$

Finally, if $v$ is proportional to a probability measure, then for any $h > 0$,

$$\mathbb{J}[Z_t] = \frac{1}{h} \mathbb{E}\left[\int_{r=0}^{r=+\infty} \int_{s=t}^{s=t+h} \int_{z=0}^{z=+\infty} N_{\mathrm{Pr}}(\mathrm{d}s, \mathrm{d}z, \mathrm{d}r)\right]$$
$$= \frac{1}{h} \times h \int_{r=0}^{r=+\infty} \int_{s=t}^{s=t+h} \int_{z=0}^{z=+\infty} r\mathbb{E}\left[X_s(r, b\pi(\mathrm{d}r))\right] v(\mathrm{d}z)$$
$$= \frac{1}{h} \times h \int_{r=0}^{r=+\infty} \int_{z=0}^{z=+\infty} r\mathbb{E}\left[X_t(r, b\pi(\mathrm{d}r))\right] v(\mathrm{d}z) \qquad . \qquad (51)$$
$$= \int_{r=0}^{r=+\infty} \int_{z=0}^{z=+\infty} r \frac{b}{r(1-M_1)} \pi(\mathrm{d}r) v(\mathrm{d}z)$$
$$= \frac{M_0}{1-M_1} b$$

**A3. Statistics of the MF model**

The strategy used here is the same as that described in **Section A2**. The statistics (14) are explicit, as in the previous model. For the average $\mathbb{E}[Z_t]$, we have

$$\mathbb{E}[Z_t] = \mathbb{E}\left[\int_{r=0}^{r=+\infty}\left\{\frac{b\pi(\mathrm{d}r)}{r}+\int_{s=-\infty}^{s=t}\int_{z=0}^{z=+\infty}e^{-r(t-s)}zN_{\mathrm{MF}}(\mathrm{d}s,\mathrm{d}z,\mathrm{d}u)\right\}\right]$$

$$= bR + M_1 \int_{r=0}^{r=+\infty}\int_{-\infty}^{t}e^{-r(t-s)}\mathbb{E}\left[wrX_{s-}(r,b\pi(\mathrm{d}r))+(1-w)\frac{b}{1-M_1}\pi(\mathrm{d}r)\right]\mathrm{d}s$$

$$= bR + M_1 \int_{r=0}^{r=+\infty}\int_{-\infty}^{t}e^{-r(t-s)}\left\{wr\frac{b}{1-M_1}\frac{1}{r}+(1-w)\frac{b}{1-M_1}\right\}\pi(\mathrm{d}r)\mathrm{d}s \qquad (52)$$

$$= bR + M_1 \frac{b}{1-M_1}\int_{r=0}^{r=+\infty}\int_{-\infty}^{t}e^{-r(t-s)}\pi(\mathrm{d}r)\mathrm{d}s$$

$$= bR + M_1 \frac{b}{1-M_1}R$$

$$= \frac{b}{1-M_1}R$$

The variance and autocovariance can be obtained again by exploiting independence among each $X(r,b\pi(\mathrm{d}r))$.

### A4. Proof

***Proof of Proposition 1***

We regard (11) as a self-consistency equation to determine $C$ as follows: for any $t \in \mathbb{R}$,

$$Cb = \mathbb{E}\left[\int_{y=0}^{y=+\infty}yX_t(y,b\pi(\mathrm{d}y))\right]$$
$$= \mathbb{E}\left[\int_{y=0}^{y=+\infty}y\left\{\frac{b}{y}\pi(\mathrm{d}y)+\int_{s=-\infty}^{s=t}\int_{z=0}^{z=+\infty}e^{-y(t-s)}zN_{\mathrm{MF}}(\mathrm{d}s,\mathrm{d}z,\mathrm{d}y)\right\}\right]. \qquad (53)$$
$$= b + \mathbb{E}\left[\int_{y=0}^{y=+\infty}y\int_{s=-\infty}^{s=t}\int_{z=0}^{z=+\infty}e^{-y(t-s)}zN_{\mathrm{MF}}(\mathrm{d}s,\mathrm{d}z,\mathrm{d}y)\right]$$

By assuming a stationary state, the right-hand side of (53) is calculated as

$$b + \mathbb{E}\left[\int_{y=0}^{y=+\infty}y\int_{s=-\infty}^{s=t}\int_{z=0}^{z=+\infty}e^{-y(t-s)}zN_{\mathrm{MF}}(\mathrm{d}s,\mathrm{d}z,\mathrm{d}y)\right]$$
$$= b + M_1\mathbb{E}\left[\int_{y=0}^{y=+\infty}y\int_{s=-\infty}^{s=t}e^{-y(t-s)}\left(ywX_s(y,b\pi(\mathrm{d}y))+(1-w)Cb\pi(\mathrm{d}y)\right)\mathrm{d}s\right]$$
$$= b + wM_1\int_{y=0}^{y=+\infty}\int_{s=-\infty}^{s=t}e^{-y(t-s)}y\mathbb{E}\left[yX_s(y,b\pi(\mathrm{d}y))\right]\mathrm{d}s + (1-w)M_1Cb\int_{y=0}^{y=+\infty}\int_{s=-\infty}^{s=t}e^{-y(t-s)}y\pi(\mathrm{d}y)\mathrm{d}s$$
$$= b + wM_1\int_{y=0}^{y=+\infty}\int_{s=-\infty}^{s=t}e^{-y(t-s)}y\mathbb{E}\left[yX_s(y,b\pi(\mathrm{d}y))\right]\mathrm{d}s + (1-w)M_1Cb$$
$$= b + wM_1\int_{y=0}^{y=+\infty}\int_{s=-\infty}^{s=t}e^{-y(t-s)}y\mathbb{E}\left[yX_t(y,b\pi(\mathrm{d}y))\right]\mathrm{d}s + (1-w)M_1Cb \leftarrow (\text{Stationarity of } X) \qquad (54)$$
$$= b + wM_1\int_{y=0}^{y=+\infty}\frac{1}{y}y\mathbb{E}\left[yX_t(y,b\pi(\mathrm{d}y))\right] + (1-w)M_1Cb$$
$$= b + wM_1\mathbb{E}\left[\int_{y=0}^{y=+\infty}yX_t(y,b\pi(\mathrm{d}y))\right] + (1-w)M_1Cb$$
$$= b + wM_1Cb + (1-w)M_1Cb$$
$$= b(1+M_1C)$$

Plugging (54) into (53) yields (13).

*Proof of Proposition 2*

The linear equation (26) can be seen as an infinite-dimensional Lyapunov equation in the context of Abi Jaber et al. [90]. We apply Theorem 3.1 of Abi Jaber et al. [90] as follows with minor modifications; note that the literature deals with a terminal-value problem while we consider an initial-value problem.

For each $\Phi \in \mathbb{L}_T^1$, define $\mathbb{T}: t \to \mathbb{T}\Phi_t$ by

$$\mathbb{T}\Phi_t = e^{-(1-wM_1)rt} + (1-w)M_1 r \int_{s=0}^{s=t} e^{-(1-wM_1)r(t-s)} \left( \int_{y=0}^{y=+\infty} \Phi_s(y) \pi(\mathrm{d}y) \right) \mathrm{d}s \text{ for } r > 0. \tag{55}$$

The equation (26) is rewritten as $E_t(1,r) = \mathbb{T}E_t(1,r)$. For any $\Phi, \Psi \in \mathbb{L}_T^1$, we have the boundedness:

$$\begin{aligned}
\|\mathbb{T}\Phi\|_{1,T} &= \sup_{0 \le t \le T} \left\| e^{-(1-wM_1)rt} + (1-w)M_1 r \int_{s=0}^{s=t} e^{-(1-wM_1)r(t-s)} \int_{y=0}^{y=+\infty} \Phi_s(y) \pi(\mathrm{d}y)\mathrm{d}s \right\|_1 \\
&\le 1 + (1-w)M_1 \sup_{0 \le t \le T} \left( \int_{r=0}^{r=+\infty} r \int_{s=0}^{s=t} e^{-(1-wM_1)r(t-s)} \|\Phi\|_{1,T} \pi(\mathrm{d}r)\mathrm{d}s \right) \\
&= 1 + (1-w)M_1 \sup_{0 \le t \le T} \left( \int_{r=0}^{r=+\infty} r \int_{s=0}^{s=t} e^{-(1-wM_1)r(t-s)} \pi(\mathrm{d}r)\mathrm{d}s \right) \|\Phi\|_{1,T} \\
&\le 1 + \frac{(1-w)M_1}{1-wM_1} \|\Phi\|_{1,T}
\end{aligned} \tag{56}$$

and hence $\|\mathbb{T}\Phi\|_{1,T,\lambda} < +\infty$ for any $\lambda > 0$. Moreover, the constant multiplied by $\|\Phi\|_{1,T}$ does not depend on $T > 0$. Similarly, the contraction follows by choosing a sufficiently large $\lambda > 0$:

$$\begin{aligned}
&\|\mathbb{T}\Phi - \mathbb{T}\Psi\|_{1,T,\lambda} \\
&= (1-w)M_1 \left\| r \int_{s=0}^{s=t} e^{-(1-wM_1)r(t-s)} \int_{y=0}^{y=+\infty} (\Phi_s(y) - \Psi_s(y)) \pi(\mathrm{d}y)\mathrm{d}s \right\|_{1,T,\lambda} \\
&\le (1-w)M_1 \sup_{0 \le t \le T} \left( e^{-\lambda t} \int_{r=0}^{r=+\infty} r \int_{s=0}^{s=t} e^{-(1-wM_1)r(t-s)} \|\Phi_s - \Psi_s\|_1 \pi(\mathrm{d}r)\mathrm{d}s \right) \\
&= (1-w)M_1 \sup_{0 \le t \le T} \left( \int_{r=0}^{r=+\infty} r \int_{s=0}^{s=t} e^{-((1-wM_1)r+\lambda)(t-s)} \left( e^{-\lambda s} \|\Phi_s - \Psi_s\|_1 \right) \pi(\mathrm{d}r)\mathrm{d}s \right). \\
&\le (1-w)M_1 \sup_{0 \le t \le T} \left( \int_{r=0}^{r=+\infty} r \int_{s=0}^{s=t} e^{-((1-wM_1)r+\lambda)(t-s)} \|\mathbb{T}\Phi - \mathbb{T}\Psi\|_{1,T,\lambda} \pi(\mathrm{d}r)\mathrm{d}s \right) \\
&\le (1-w)M_1 \|\mathbb{T}\Phi - \mathbb{T}\Psi\|_{1,T,\lambda} \sup_{0 \le t \le T} \left( \int_{r=0}^{r=+\infty} r \int_{s=0}^{s=t} e^{-((1-wM_1)r+\lambda)(t-s)} \pi(\mathrm{d}r)\mathrm{d}s \right) \\
&\le \left\{ (1-w)M_1 \int_{r=0}^{r=+\infty} \frac{r}{(1-wM_1)r+\lambda} \pi(\mathrm{d}r) \right\} \|\mathbb{T}\Phi - \mathbb{T}\Psi\|_{1,T,\lambda}
\end{aligned} \tag{57}$$

By the classical dominated convergence theorem, the coefficient inside $\{\cdot\}$ in the last line of (57) can be made strictly smaller than 1 by choosing a sufficiently large $\lambda > 0$. Moreover, this coefficient does not depend on $T > 0$. Therefore, the equation $E_t(1,\cdot) = \mathbb{T}E_t(1,\cdot)$ ($t \ge 0$) admits a unique solution belonging to the space $\mathbb{L}_T^1$ for any $T > 0$. The facts that this solution satisfies the Lyapunov equation (26) and the initial condition $E_0(1,\cdot) = 1$ can be checked by a direct substitution. Along with (58) below, the right-hand side of (55) is continuous for $t > 0$.

The well-posedness of the second Lyapunov equation (27) follows in the same way because of

$$\|E(1,\cdot)\|_{1,T} \leq \frac{1-wM_1}{1-M_1} \quad (58)$$

for any $T > 0$ which follows from $\|E(1,\cdot)\|_{1,T} = \|\mathbb{T}E(1,\cdot)\|_{1,T}$ and (56), and the pointwise bound due to (58):

$$\begin{aligned}
|E_t(1,r)| &= \left| e^{-(1-wM_1)rt} + (1-w)M_1 r \int_{s=0}^{s=t} e^{-(1-wM_1)r(t-s)} \left( \int_{y=0}^{y=+\infty} E_s(1,y)\pi(\mathrm{d}y) \right) \mathrm{d}s \right| \\
&\leq 1 + (1-w)M_1 r \int_{s=0}^{s=t} e^{-(1-wM_1)r(t-s)} \left( \int_{y=0}^{y=+\infty} |E_s(1,y)|\pi(\mathrm{d}y) \right) \mathrm{d}s \\
&\leq 1 + (1-w)M_1 r \int_{s=0}^{s=t} e^{-(1-wM_1)r(t-s)} \mathrm{d}s \|E(1,\cdot)\|_{1,T} \\
&\leq 1 + (1-w)M_1 \times \frac{1-wM_1}{1-M_1} r \int_{s=0}^{s=t} e^{-(1-wM_1)r(t-s)} \mathrm{d}s \\
&\leq 1 + (1-w)M_1 \times \frac{1-wM_1}{1-M_1} r \times \frac{1}{(1-wM_1)r} \\
&= 1 + \frac{(1-w)M_1}{1-M_1} \\
&= \frac{1-wM_1}{1-M_1}
\end{aligned} \quad (59)$$

for any $T \geq t \geq 0$ and $r > 0$.

□

*Proof of Proposition 3*

The proof follows the lines of **Proof of Proposition 2** with technical modifications to deal with the nonlinearity of the generalized Riccati equation. For this purpose, we start from an auxiliary equation with a truncated nonlinear term. Below, $L = \frac{\theta}{1-M_1} > 0$ is a constant. For each $\Phi \in \mathbb{L}_T^1$ and $T > 0$, set $|\Phi|_{\infty,T} = \sup \Phi_t(r)$ where the supremum is taken for $r > 0$ and $0 \leq t \leq T$.

For each $\Phi \in \mathbb{L}_T^1$ and $T > 0$, fix one $\hat{\Phi} \in \mathbb{L}_T^1$ such that $|\hat{\Phi}|_{\infty,T} \leq L$, and define $\mathbb{U}: t \to \mathbb{U}\Phi_t$ by

$$\mathbb{U}\Phi_t = \theta e^{-rt} + r \int_{s=0}^{s=t} e^{-r(t-s)} \left\{ \begin{array}{l} w \int_0^{+\infty} \left( 1 - e^{-(\hat{\Phi}_s(r))_+ z} \right) \nu(\mathrm{d}z) \\ + (1-w) \int_{y=0}^{y=+\infty} \int_{z=0}^{z=+\infty} \left( 1 - e^{-(\Phi_s(y))_+ z} \right) \nu(\mathrm{d}z)\pi(\mathrm{d}y) \end{array} \right\} \mathrm{d}s \text{ for } r > 0. \quad (60)$$

Here, $(x)_+ = \max\{x, 0\}$ for any $x \in \mathbb{R}$. This equation comes from the modified Riccati equation

$$\frac{\partial}{\partial t} B_t(r,\theta) = r \left\{ \begin{array}{l} -B_t(r,\theta) + w \int_0^{+\infty} \left( 1 - e^{-(B_t(r,\theta))_+ z} \right) \nu(\mathrm{d}z) \\ + (1-w) \int_{y=0}^{y=+\infty} \int_{z=0}^{z=+\infty} \left( 1 - e^{-(B_t(y,\theta))_+ z} \right) \nu(\mathrm{d}z)\pi(\mathrm{d}y) \end{array} \right\}, \; r > 0 \text{ and } t > 0 \quad (61)$$

subject to the initial condition $B_0(\theta, \cdot) = \theta$, because we can formally rewrite it as $B_t = \mathbb{U}B_t$ if $B$ is not larger than $L$ pointwise and $\Phi = B$. The generalized Riccati equation can be expressed similarly by using $\mathbb{U}$ with dropping "max". A strategy here is to show that the modified Riccati equation admits a unique solution in $\mathbb{L}_T^1$ for any $T > 0$, and this solution is nonnegative and bounded above by $L$ pointwise. Then, we show that this solution is the desired solution to the generalized Riccati equation (22), and the taking "max" in (61) turns out to be actually superficial.

As in **Proof of Proposition 2**, we show the boundedness and contraction property of $\mathbb{U}$. They are proven as follows. For each $T > 0$, temporally fix one $\hat{\Phi} \in \mathbb{L}_T^1$ such that $\left|\hat{\Phi}\right|_{\infty,T} \leq L$. For each $\Phi \in \mathbb{L}_T^1$, by $1 - e^{-ax} \leq ax$ for any $a \geq 0$ and $x \geq 0$, we have

$$
\begin{aligned}
\|\mathbb{U}\Phi\|_{1,T} &= \sup_{0 \leq t \leq T} \left\| \theta e^{-rt} + r\int_{s=0}^{s=t} e^{-r(t-s)} \left\{ \begin{array}{l} w\int_0^{+\infty}\left(1 - e^{-(\hat{\Phi}_s(r))z}\right)\nu(\mathrm{d}z) \\ + (1-w)\int_{y=0}^{y=+\infty}\int_{z=0}^{z=+\infty}\left(1 - e^{-(\Phi_s(y))_+ z}\right)\nu(\mathrm{d}z)\pi(\mathrm{d}y) \end{array} \right\} \mathrm{d}s \right\|_1 \\
&\leq \theta + \sup_{0 \leq t \leq T}\left\| r\int_{s=0}^{s=t} e^{-r(t-s)} \left\{ \begin{array}{l} w\int_0^{+\infty}\left(1 - e^{-(\hat{\Phi}_s(r))z}\right)\nu(\mathrm{d}z) \\ + (1-w)\int_{y=0}^{y=+\infty}\int_{z=0}^{z=+\infty}\left(1 - e^{-(\Phi_s(y))_+ z}\right)\nu(\mathrm{d}z)\pi(\mathrm{d}y) \end{array} \right\} \mathrm{d}s \right\|_1 \\
&\leq \theta + \sup_{0 \leq t \leq T}\left\| r\int_{s=0}^{s=t} e^{-r(t-s)} \left\{ \begin{array}{l} w\left(\hat{\Phi}_s(r)\right)\int_0^{+\infty} z\nu(\mathrm{d}z) \\ + (1-w)\left(\int_{y=0}^{y=+\infty}(\Phi_s(y))_+\pi(\mathrm{d}y)\right)\int_{z=0}^{z=+\infty} z\nu(\mathrm{d}z) \end{array} \right\} \mathrm{d}s \right\|_1 \\
&\leq \theta + M_1 \sup_{0 \leq t \leq T}\left\| r\int_{s=0}^{s=t} e^{-r(t-s)}\left\{ wL + (1-w)\|\Phi\|_{1,T} \right\} \mathrm{d}s \right\| \\
&= \theta + \left( wL + (1-w)\|\Phi\|_{1,T} \right) M_1 \sup_{0 \leq t \leq T}\left\| r\int_{s=0}^{s=t} e^{-r(t-s)} \mathrm{d}s \right\| \\
&\leq \theta + \left( wL + (1-w)\|\Phi\|_{1,T} \right) M_1
\end{aligned}
\tag{62}
$$

and hence $\|\mathbb{U}\Phi\|_{1,T,\lambda} < +\infty$ for any $\lambda > 0$ independent from $T$. Similarly, for any $\Phi, \Psi \in \mathbb{L}_T^1$, the contraction follows by choosing a sufficiently large $\lambda > 0$ (we use $\left|e^{-ax} - e^{-ay}\right| \leq a|x - y|$ for any $a \geq 0$ and $x, y \geq 0$):

$$
\begin{aligned}
&\|\mathbb{U}\Phi - \mathbb{U}\Psi\|_{1,T,\lambda} \\
&= (1-w)\sup_{0 \leq t \leq T}\left\| r\int_{s=0}^{s=t} e^{-r(t-s)} \left\{ \int_{y=0}^{y=+\infty}\int_{z=0}^{z=+\infty}\left(e^{-(\Psi_s(y))_+ z} - e^{-(\Phi_s(y))_+ z}\right)\nu(\mathrm{d}z)\pi(\mathrm{d}y) \right\} \mathrm{d}s \right\|_{1,T,\lambda} \\
&\leq (1-w)\sup_{0 \leq t \leq T}\left\| r\int_{s=0}^{s=t} e^{-r(t-s)} \left\{ \left(\int_{y=0}^{y=+\infty}\left|(\Psi_s(y))_+ - (\Phi_s(y))_+\right|\pi(\mathrm{d}y)\right)\int_{z=0}^{z=+\infty} z\nu(\mathrm{d}z) \right\} \mathrm{d}s \right\|_{1,T,\lambda} \\
&\leq (1-w)M_1 \sup_{0 \leq t \leq T}\left\| r\int_{s=0}^{s=t} e^{-r(t-s)}\|\Phi_s - \Psi_s\|_1 \mathrm{d}s \right\|_{1,T,\lambda} \\
&= (1-w)M_1 \sup_{0 \leq t \leq T}\left( e^{-\lambda t}\int_{r=0}^{r=+\infty} r\int_{s=0}^{s=t} e^{-r(t-s)}\|\Phi_s - \Psi_s\|_1 \mathrm{d}s \pi(\mathrm{d}r) \right) \\
&= (1-w)M_1 \sup_{0 \leq t \leq T}\left( \int_{r=0}^{r=+\infty} r\int_{s=0}^{s=t} e^{-(r+\lambda)(t-s)}\left(e^{-\lambda s}\|\Phi_s - \Psi_s\|_1\right) \mathrm{d}s \pi(\mathrm{d}r) \right) \\
&\leq (1-w)M_1 \|\Phi - \Psi\|_{1,T,\lambda}\sup_{0 \leq t \leq T}\left( \int_{r=0}^{r=+\infty} r\int_{s=0}^{s=t} e^{-(r+\lambda)(t-s)}\mathrm{d}s\pi(\mathrm{d}r) \right) \\
&\leq \left\{(1-w)M_1 \int_0^{+\infty}\frac{r}{r+\lambda}\pi(\mathrm{d}r)\right\}\|\Phi - \Psi\|_{1,T,\lambda}
\end{aligned}
\tag{63}
$$

where the coefficient multiplied by $\|\Phi-\Psi\|_{1,T,\lambda}$ becomes strictly smaller than 1 by choosing some $\lambda>0$. Consequently, as in **Proof of Proposition 2**, the boundedness and contraction show that the following equation admits a unique solution in $\mathbb{L}_T^1$ such that $|\Phi|_{\infty,T}\le L$:

$$\Phi_t(r)=\mathbb{U}\Phi_t(r),\ 0\le t\le T\ \text{and}\ r>0. \tag{64}$$

Moreover, this solution is bounded pointwise because of

$$|\Phi|_{\infty,T}=\left|\theta e^{-rt}+r\int_{s=0}^{s=t}e^{-r(t-s)}\left\{\begin{array}{l}w\int_0^{+\infty}\left(1-e^{-(\hat{\Phi}_s(r))z}\right)\nu(\mathrm{d}z)\\+(1-w)\int_{y=0}^{y=+\infty}\int_{z=0}^{z=+\infty}\left(1-e^{-(\Phi_s(y))_+z}\right)\nu(\mathrm{d}z)\pi(\mathrm{d}y)\end{array}\right\}\mathrm{d}s\right|_{\infty,T}$$

$$\le \theta+w\left|r\int_{s=0}^{s=t}e^{-r(t-s)}\int_0^{+\infty}\left(1-e^{-(\hat{\Phi}_s(r))z}\right)\nu(\mathrm{d}z)\mathrm{d}s\right|_{\infty,T} \tag{65}$$

$$+(1-w)\left|r\int_{s=0}^{s=t}e^{-r(t-s)}\int_{y=0}^{y=+\infty}\int_{z=0}^{z=+\infty}\left(1-e^{-(\Phi_s(y))_+z}\right)\nu(\mathrm{d}z)\pi(\mathrm{d}y)\mathrm{d}s\right|_{\infty,T}$$

$$\le \theta+wLM_1+(1-w)M_1|\Phi|_{\infty,T}$$

and hence

$$|\Phi|_{\infty,T}\le\frac{\theta+wLM_1}{1-(1-w)M_1}=\frac{\theta}{1-M_1}. \tag{66}$$

Similarly, by (64) we have

$$\|\Phi\|_{1,T}\le\frac{\theta}{1-M_1}. \tag{67}$$

Fix $T>0$. We consider a sequence of mappings $\Phi^{(n)}\in\mathbb{L}_T^1$ ($n=0,1,2,3,\ldots$) with a constant function $\Phi^{(0)}\equiv\rho\in[0,L]$ to apply a fixed-point argument. The operator $\mathbb{U}$ with $\hat{\Phi}=\Phi^{(n)}$ is expressed as $\mathbb{U}^{(n)}$, and we consider the following recursion for $n=0,1,2,\ldots$:

$$\Phi_t^{(n)}(r)=\mathbb{U}^{(n-1)}\Phi_t^{(n)}(r),\ 0\le t\le T\ \text{and}\ r>0. \tag{68}$$

By the contraction and boundedness property of $\mathbb{U}^{(0)}$, there exists a unique $\Phi^{(1)}\in\mathbb{L}_T^1$ satisfying the uniform bounds (66)-(67). We can continue this procedure for larger $n$ in a cascading manner, and obtain a unique sequence $\Phi^{(n)}\in\mathbb{L}_T^1$ ($n=0,1,2,\ldots$) all elements satisfying the uniform bounds (66)-(67). Then, applying a Banach's fixed point theorem (e.g., Section 4 of Abi Jaber et al. [90]) shows that there exists a unique solution $\Phi\in\mathbb{L}_T^1$ to

$$\Phi_t(r)=\overline{\mathbb{U}}\Phi_t(r)=\theta e^{-rt}+r\int_{s=0}^{s=t}e^{-r(t-s)}\left\{\begin{array}{l}w\int_0^{+\infty}\left(1-e^{-(\Phi_s(r))_+z}\right)\nu(\mathrm{d}z)\\+(1-w)\int_{y=0}^{y=+\infty}\int_{z=0}^{z=+\infty}\left(1-e^{-(\Phi_s(y))_+z}\right)\nu(\mathrm{d}z)\pi(\mathrm{d}y)\end{array}\right\}\mathrm{d}s \tag{69}$$

for $t\ge 0$ ($T$ is arbitrary) and $r>0$, satisfying the uniform bounds (66)-(67). This solution, denoted by $\Phi$, is nonnegative. Indeed, for any $t\ge 0$ and $r>0$, by the elementary inequality $1-e^{-x}\ge 0$ for any $x\ge 0$,

$$\Phi_t(r)=\overline{\mathbb{U}}\Phi_t(r)\ge\theta e^{-rt}\ge 0. \tag{70}$$

The inequality (70) shows that the solution to the modified Riccati equation (61), which is nonnegative, also satisfies the generalized Riccati equation (22). This is the desired unique nonnegative solution by the contraction property. Finally, on the continuity, for the solution $\Phi$ found above, for any $t, u \geq 0$ we have

$$\left\| \overline{\mathbb{U}} \Phi_t - \overline{\mathbb{U}} \Phi_u \right\|_1$$

$$= \left\| \begin{array}{l} r \int_{s=0}^{s=t} e^{-r(t-s)} \left\{ w \int_0^{+\infty} \left(1 - e^{-(\Phi_s(r))_+ z}\right) v(\mathrm{d}z) \\ + (1-w) \int_{y=0}^{y=+\infty} \int_{z=0}^{z=+\infty} \left(1 - e^{-(\Phi_s(y))_+ z}\right) v(\mathrm{d}z) \pi(\mathrm{d}y) \right\} \mathrm{d}s \\ -r \int_{s=0}^{s=u} e^{-r(u-s)} \left\{ w \int_0^{+\infty} \left(1 - e^{-(\Phi_s(r))_+ z}\right) v(\mathrm{d}z) \\ + (1-w) \int_{y=0}^{y=+\infty} \int_{z=0}^{z=+\infty} \left(1 - e^{-(\Phi_s(y))_+ z}\right) v(\mathrm{d}z) \pi(\mathrm{d}y) \right\} \mathrm{d}s \end{array} \right\|_1 \quad (71)$$

The right-hand side of (71) vanishes as $u \to t$ due to a uniform boundedness of the quantities inside "$\{\cdot\}$" and dominated convergence that follows from (66)-(67) and

$$\left\| r \int_0^t e^{-r(t-s)} \mathrm{d}s \right\|, \left\| r \int_0^u e^{-r(u-s)} \mathrm{d}s \right\| \leq \left\| r \times \frac{1}{r} \right\| = 1 . \quad (72)$$

□

***Proof of (28)***

First, the left equation of (28) is obtained as follows. By (26), we have

$$\int_{t=0}^{t=+\infty} \int_{r=0}^{r=+\infty} \frac{1}{r} \frac{\partial}{\partial t} E_t(1, r) \pi(\mathrm{d}r) \mathrm{d}t$$
$$= \int_{t=0}^{t=+\infty} \int_{r=0}^{r=+\infty} \left\{ -(1 - wM_1) E_t(1, r) + (1-w) M_1 \int_{y=0}^{y=+\infty} E_t(1, y) \pi(\mathrm{d}y) \right\} \pi(\mathrm{d}r) \mathrm{d}t \quad (73)$$

If $E_t(1, \cdot) \to 0$ as $t \to +\infty$, then its left-hand side becomes (by exchanging the order of integrals)

$$\int_{t=0}^{t=+\infty} \int_{r=0}^{r=+\infty} \frac{1}{r} \frac{\partial}{\partial t} E_t(1, r) \pi(\mathrm{d}r) \mathrm{d}t = \int_{r=0}^{r=+\infty} \frac{1}{r} \left( \int_{t=0}^{t=+\infty} \frac{\partial}{\partial t} E_t(1, r) \mathrm{d}t \right) \pi(\mathrm{d}r)$$
$$= \int_{r=0}^{r=+\infty} \frac{1}{r} \left[ E_t(1, r) \right]_{t=0}^{t=+\infty} \pi(\mathrm{d}r) \quad (74)$$
$$= -\int_{r=0}^{r=+\infty} \frac{1}{r} \pi(\mathrm{d}r)$$
$$= -R$$

The right-hand side is rewritten as

$$\int_{t=0}^{t=+\infty} \int_{r=0}^{r=+\infty} \left\{ -(1 - wM_1) E_t(1, r) + (1-w) M_1 \int_{y=0}^{y=+\infty} E_t(1, y) \pi(\mathrm{d}y) \right\} \pi(\mathrm{d}r) \mathrm{d}t$$
$$= -(1 - wM_1) \int_{t=0}^{t=+\infty} \int_{r=0}^{r=+\infty} E_t(1, r) \pi(\mathrm{d}r) \mathrm{d}t + (1-w) M_1 \int_{t=0}^{t=+\infty} \int_{r=0}^{r=+\infty} \int_{y=0}^{y=+\infty} E_t(1, y) \pi(\mathrm{d}y) \pi(\mathrm{d}r) \mathrm{d}t \quad (75)$$
$$= \{-(1 - wM_1) + (1-w) M_1\} \int_{t=0}^{t=+\infty} \int_{r=0}^{r=+\infty} E_t(1, r) \pi(\mathrm{d}r) \mathrm{d}t$$
$$= -(1 - M_1) \int_{t=0}^{t=+\infty} \int_{r=0}^{r=+\infty} E_t(1, r) \pi(\mathrm{d}r) \mathrm{d}t$$

Substituting (74) and (75) into (73) yields the desired result.

Second, the right equation of (28) is proven as follows. By (27), we have

$$\int_{t=0}^{t=+\infty}\int_{r=0}^{r=+\infty} \frac{1}{r}\frac{\partial}{\partial t} E_t(2,r) \pi(\mathrm{d}r)\mathrm{d}t$$
$$= \int_{t=0}^{t=+\infty}\int_{r=0}^{r=+\infty} \left[ \begin{array}{l} -(1-wM_1)E_t(2,r)+(1-w)M_1 \int_{y=0}^{y=+\infty} E_t(2,y)\pi(\mathrm{d}y) \\ +M_2\left\{w(E_t(1,r))^2+(1-w)\int_{y=0}^{y=+\infty}(E_t(1,y))^2 \pi(\mathrm{d}y)\right\} \end{array}\right] \pi(\mathrm{d}r)\mathrm{d}t \quad (76)$$

If $E_t(2,\cdot) \to 0$ as $t \to +\infty$, then its left-hand side is rewritten as (by exchanging the order of integrals)

$$\int_{t=0}^{t=+\infty}\int_{r=0}^{r=+\infty} \frac{1}{r}\frac{\partial}{\partial t} E_t(2,r) \pi(\mathrm{d}r)\mathrm{d}t = 0. \quad (77)$$

The right-hand side of (76) becomes

$$\int_{t=0}^{t=+\infty}\int_{r=0}^{r=+\infty} \left[ \begin{array}{l} -(1-wM_1)E_t(2,r)+(1-w)M_1 \int_{y=0}^{y=+\infty} E_t(2,y)\pi(\mathrm{d}y) \\ +M_2\left\{w(E_t(1,r))^2+(1-w)\int_{y=0}^{y=+\infty}(E_t(1,y))^2 \pi(\mathrm{d}y)\right\} \end{array}\right] \pi(\mathrm{d}r)\mathrm{d}t$$
$$= -(1-M_1)\int_{t=0}^{t=+\infty}\int_{r=0}^{r=+\infty} E_t(2,r)\pi(\mathrm{d}r)\mathrm{d}t \quad (78)$$
$$+M_2 \int_{t=0}^{t=+\infty}\int_{r=0}^{r=+\infty} \left\{w(E_t(1,r))^2+(1-w)\int_{y=0}^{y=+\infty}(E_t(1,y))^2 \pi(\mathrm{d}y)\right\}\pi(\mathrm{d}r)\mathrm{d}t$$
$$= -(1-M_1)\int_{t=0}^{t=+\infty}\int_{r=0}^{r=+\infty} E_t(2,r)\pi(\mathrm{d}r)\mathrm{d}t + M_2 \int_{t=0}^{t=+\infty}\int_{r=0}^{r=+\infty}(E_t(1,r))^2 \pi(\mathrm{d}r)\mathrm{d}t$$

Substituting (77) and (78) into (76) yields

$$\mathbb{V}[Z_t] = \frac{bM_2}{1-M_1}\int_{t=0}^{t=+\infty}\int_{r=0}^{r=+\infty}(E_t(1,r))^2 \pi(\mathrm{d}r)\mathrm{d}t. \quad (79)$$

We need another step to yield the desired result. We have

$$\frac{1}{2}\int_{t=0}^{t=+\infty}\int_{r=0}^{r=+\infty} \frac{1}{r}\frac{\partial}{\partial t}(E_t(1,r))^2 \pi(\mathrm{d}r)\mathrm{d}t$$
$$= \int_{t=0}^{t=+\infty}\int_{r=0}^{r=+\infty}\left\{-(1-wM_1)(E_t(1,r))^2+(1-w)M_1 E_t(1,r)\int_{y=0}^{y=+\infty}E_t(1,y)\pi(\mathrm{d}y)\right\}\pi(\mathrm{d}r)\mathrm{d}t \quad (80)$$

and hence

$$-\frac{1}{2}R = -(1-wM_1)\int_{t=0}^{t=+\infty}\int_{r=0}^{r=+\infty}(E_t(1,r))^2 \pi(\mathrm{d}r)\mathrm{d}t+(1-w)M_1\int_{t=0}^{t=+\infty}\left(\int_{r=0}^{r=+\infty}E_t(1,r)\pi(\mathrm{d}r)\right)^2 \mathrm{d}t. \quad (81)$$

We then obtain

$$\int_{t=0}^{t=+\infty}\int_{r=0}^{r=+\infty}(E_t(1,r))^2 \pi(\mathrm{d}r)\mathrm{d}t = \frac{1}{1-wM_1}\left\{\frac{1}{2}R+(1-w)M_1\int_{t=0}^{t=+\infty}\left(\int_{r=0}^{r=+\infty}E_t(1,r)\pi(\mathrm{d}r)\right)^2 \mathrm{d}t\right\}. \quad (82)$$

Substituting (82) into (79) completes the proof.

# Supplementary file of "Superposition of interacting stochastic processes with memory and its application to migrating fish counts" by Hidekazu Yoshioka

**A5. Biological count data used in the main text**

The daily count data of *P. altivelis* are summarized in this section. **Table S1** reports the total number of analyzed data points $I$, and total number of fish counts for each year and river. Second, **Table S2** summarizes the empirical average (Ave), variance (Var), coefficient of variation (CV), total number of jumps (Jmp), and skewness (Skw) of each time series data. **Table S3** presents identified parameter values, with $w=1$ being fixed. **Table S4** details the parameter values identified with $w$ being fitted. **Table S5** presents empirical skewness with $w$ being fixed and that with $w$ being fitted.

We investigated the correlations between the $J$ of the Yahagi, Tama, Ara, and Tone Rivers between 2012 and 2023, as reported in **Table S6**; however, it is not the focus of this study. The correlations between the Tama and Ara Rivers ($R^2$ value of 0.9272) and that between the Tama and Tone Rivers ($R^2$ value of 0.4268) are higher than the others and are smaller than 0.1. Therefore, the total number of migrants is not always correlated among these rivers; however, they all enter the Pacific Ocean, and three of them (Tama, Ara, and Tone Rivers) flow in the Kanto Region in Japan. The high correlation between the Tama and Ara Rivers is attributed to their flow into the same bay.

Table S1. The total number of data points $I$ and fish counts $J$ for each year and river. **Table 1** presents the data sources, and the author determined the total number of fish counts $J$. "First day" in the table represents the first day with a positive fish count during an observation period.

| River | Year | First day | $I$ | $J$ |
|---|---|---|---|---|
| Yahagi River | 1998 | April 20 | 42 | 3112568 |
| | 1999 | April 15 | 62 | 223632 |
| | 2000 | May 5 | 35 | 53317 |
| | 2001 | May 23 | 25 | 1289 |
| | 2002 | May 5 | 26 | 218713 |
| | 2003 | April 21 | 78 | 316020 |
| | 2004 | April 1 | 91 | 2121313 |
| | 2005 | April 6 | 100 | 559637 |
| | 2006 | April 27 | 65 | 712433 |
| | 2007 | March 29 | 95 | 6218898 |
| | 2008 | March 24 | 101 | 637088 |
| | 2009 | April 8 | 65 | 1081435 |
| | 2010 | April 16 | 91 | 487951 |
| | 2011 | April 10 | 97 | 985637 |
| | 2012 | April 13 | 91 | 761990 |
| | 2013 | March 30 | 108 | 839587 |
| | 2014 | April 15 | 88 | 601147 |
| | 2015 | April 12 | 89 | 1276048 |
| | 2016 | April 11 | 91 | 10030840 |
| | 2017 | April 12 | 80 | 1440609 |
| | 2018 | April 12 | 80 | 2307520 |
| | 2019 | April 17 | 54 | 447134 |
| | 2020 | April 7 | 84 | 1103486 |
| | 2021 | April 6 | 40 | 603673 |
| | 2022 | April 2 | 90 | 913896 |
| | 2023 | April 20 | 42 | 43127 |
| Nagara River | 2023 | February 22 | 127 | 852596 |
| | 2024 | February 28 | 120 | 1236102 |
| Tama River | 2011 | March 25 | 66 | 422585 |
| | 2012 | March 23 | 66 | 644779 |
| | 2013 | March 19 | 74 | 348081 |
| | 2014 | March 18 | 75 | 292075 |
| | 2015 | March 20 | 73 | 234760 |
| | 2016 | March 23 | 70 | 250193 |
| | 2017 | March 17 | 76 | 85487 |
| | 2018 | March 20 | 72 | 536528 |
| | 2019 | March 23 | 70 | 179669 |
| | 2020 | March 25 | 58 | 19923 |
| | 2021 | March 17 | 69 | 17518 |
| | 2022 | March 9 | 84 | 135117 |
| | 2023 | March 10 | 81 | 112108 |
| | 2024 | March 15 | 64 | 19880 |
| Ara River | 2012 | April 1 | 47 | 899130 |
| | 2013 | April 8 | 40 | 441279 |
| | 2014 | April 6 | 40 | 408802 |
| | 2015 | April 6 | 47 | 146830 |
| | 2016 | April 6 | 40 | 158099 |
| | 2017 | April 6 | 40 | 222545 |
| | 2018 | April 6 | 40 | 530176 |
| | 2019 | April 6 | 40 | 44768 |
| | 2020 | April 6 | 40 | 60731 |

|  | 2021 | April 6 | 40 | 44663 |
|---|---|---|---|---|
|  | 2022 | March 27 | 50 | 191328 |
|  | 2023 | March 29 | 48 | 106407 |
|  | 2024 | March 29 | 43 | 61725 |
| Tone River | 2012 | April 21 | 41 | 101471 |
|  | 2013 | April 26 | 36 | 17282 |
|  | 2014 | April 22 | 40 | 81488 |
|  | 2015 | April 27 | 35 | 10845 |
|  | 2016 | April 21 | 41 | 14874 |
|  | 2017 | April 22 | 40 | 35744 |
|  | 2018 | April 21 | 41 | 70847 |
|  | 2019 | April 21 | 41 | 15756 |
|  | 2020 | April 22 | 39 | 54683 |
|  | 2021 | April 21 | 42 | 15976 |
|  | 2022 | April 21 | 41 | 189564 |
|  | 2023 | April 21 | 41 | 85485 |
|  | 2024 | April 21 | 41 | 15303 |

Table S2. Empirical average (Ave), variance (Var), coefficient of variation (CV), total number of jumps (Jmp), and skewness (Skw) of each data.

| River | Year | Ave (ind/day) | Var (ind$^2$/day$^2$) | CV (-) | Jmp (-) | Skw (-) |
|---|---|---|---|---|---|---|
| Yahagi River | 1998 | 8.191.E+04 | 1.362.E+10 | 1.425.E+00 | 3.095.E-01 | 1.550.E+00 |
| | 1999 | 3.607.E+03 | 1.286.E+08 | 3.144.E+00 | 2.903.E-01 | 4.954.E+00 |
| | 2000 | 1.523.E+03 | 1.184.E+07 | 2.259.E+00 | 3.143.E-01 | 2.928.E+00 |
| | 2001 | 2.148.E+02 | 1.315.E+05 | 1.688.E+00 | 1.200.E-01 | 1.687.E+00 |
| | 2002 | 8.412.E+03 | 3.746.E+08 | 2.301.E+00 | 2.692.E-01 | 3.350.E+00 |
| | 2003 | 4.052.E+03 | 1.760.E+08 | 3.274.E+00 | 3.077.E-01 | 4.627.E+00 |
| | 2004 | 2.331.E+04 | 1.613.E+10 | 5.448.E+00 | 2.308.E-01 | 6.448.E+00 |
| | 2005 | 6.083.E+03 | 1.728.E+08 | 2.161.E+00 | 2.700.E-01 | 2.863.E+00 |
| | 2006 | 1.113.E+04 | 9.126.E+08 | 2.714.E+00 | 2.923.E-01 | 4.462.E+00 |
| | 2007 | 6.687.E+04 | 3.460.E+10 | 2.782.E+00 | 2.947.E-01 | 5.292.E+00 |
| | 2008 | 9.509.E+03 | 4.044.E+08 | 2.115.E+00 | 1.980.E-01 | 3.360.E+00 |
| | 2009 | 2.207.E+04 | 1.760.E+09 | 1.901.E+00 | 2.000.E-01 | 2.536.E+00 |
| | 2010 | 5.483.E+03 | 1.234.E+08 | 2.026.E+00 | 3.077.E-01 | 3.350.E+00 |
| | 2011 | 1.038.E+04 | 3.523.E+08 | 1.809.E+00 | 3.299.E-01 | 3.287.E+00 |
| | 2012 | 8.374.E+03 | 1.874.E+08 | 1.635.E+00 | 3.297.E-01 | 2.304.E+00 |
| | 2013 | 7.774.E+03 | 1.897.E+08 | 1.772.E+00 | 2.963.E-01 | 4.618.E+00 |
| | 2014 | 6.910.E+03 | 1.687.E+08 | 1.880.E+00 | 2.841.E-01 | 3.448.E+00 |
| | 2015 | 1.484.E+04 | 7.532.E+08 | 1.850.E+00 | 2.921.E-01 | 2.855.E+00 |
| | 2016 | 1.102.E+05 | 2.696.E+10 | 1.490.E+00 | 2.198.E-01 | 1.567.E+00 |
| | 2017 | 1.801.E+04 | 7.427.E+08 | 1.513.E+00 | 2.250.E-01 | 2.550.E+00 |
| | 2018 | 3.077.E+04 | 3.273.E+09 | 1.860.E+00 | 3.125.E-01 | 2.981.E+00 |
| | 2019 | 8.436.E+03 | 1.069.E+08 | 1.225.E+00 | 3.148.E-01 | 2.339.E+00 |
| | 2020 | 1.314.E+04 | 3.465.E+08 | 1.417.E+00 | 9.524.E-02 | 1.854.E+00 |
| | 2021 | 1.548.E+04 | 9.041.E+08 | 1.943.E+00 | 2.750.E-01 | 2.927.E+00 |
| | 2022 | 1.039.E+04 | 3.719.E+08 | 1.857.E+00 | 2.000.E-01 | 2.901.E+00 |
| | 2023 | 1.027.E+03 | 4.257.E+06 | 2.009.E+00 | 1.905.E-01 | 2.791.E+00 |
| Nagara River | 2023 | 7.287.E+03 | 2.136.E+08 | 2.006.E+00 | 2.756.E-01 | 2.672.E+00 |
| | 2024 | 1.189.E+04 | 5.540.E+08 | 1.980.E+00 | 2.583.E-01 | 2.713.E+00 |
| Tama River | 2011 | 6.403.E+03 | 2.573.E+08 | 2.505.E+00 | 2.879.E-01 | 5.924.E+00 |
| | 2012 | 9.920.E+03 | 2.010.E+08 | 1.429.E+00 | 3.485.E-01 | 2.239.E+00 |
| | 2013 | 4.903.E+03 | 4.059.E+07 | 1.300.E+00 | 2.838.E-01 | 1.774.E+00 |
| | 2014 | 3.947.E+03 | 2.634.E+07 | 1.300.E+00 | 2.933.E-01 | 2.056.E+00 |
| | 2015 | 3.306.E+03 | 1.950.E+07 | 1.335.E+00 | 3.425.E-01 | 1.512.E+00 |
| | 2016 | 3.626.E+03 | 3.964.E+07 | 1.736.E+00 | 2.714.E-01 | 2.889.E+00 |
| | 2017 | 1.125.E+03 | 2.789.E+06 | 1.485.E+00 | 3.026.E-01 | 3.567.E+00 |
| | 2018 | 8.254.E+03 | 3.983.E+08 | 2.418.E+00 | 2.500.E-01 | 6.434.E+00 |
| | 2019 | 2.604.E+03 | 5.959.E+07 | 2.964.E+00 | 2.714.E-01 | 5.425.E+00 |
| | 2020 | 3.985.E+02 | 6.222.E+05 | 1.980.E+00 | 2.759.E-01 | 3.514.E+00 |
| | 2021 | 2.654.E+02 | 2.319.E+05 | 1.814.E+00 | 2.899.E-01 | 3.524.E+00 |
| | 2022 | 1.851.E+03 | 6.273.E+06 | 1.353.E+00 | 2.262.E-01 | 2.021.E+00 |
| | 2023 | 1.437.E+03 | 4.688.E+06 | 1.506.E+00 | 2.840.E-01 | 2.441.E+00 |
| | 2024 | 3.488.E+02 | 2.316.E+05 | 1.380.E+00 | 2.344.E-01 | 1.479.E+00 |
| Ara River | 2012 | 1.998.E+04 | 1.076.E+09 | 1.642.E+00 | 2.766.E-01 | 3.247.E+00 |
| | 2013 | 1.103.E+04 | 2.356.E+08 | 1.391.E+00 | 3.000.E-01 | 1.942.E+00 |
| | 2014 | 1.022.E+04 | 1.771.E+08 | 1.302.E+00 | 3.000.E-01 | 1.972.E+00 |
| | 2015 | 3.124.E+03 | 9.384.E+06 | 9.806.E-01 | 3.191.E-01 | 1.144.E+00 |
| | 2016 | 3.952.E+03 | 2.129.E+07 | 1.167.E+00 | 3.500.E-01 | 1.911.E+00 |
| | 2017 | 5.564.E+03 | 3.462.E+07 | 1.058.E+00 | 3.250.E-01 | 1.773.E+00 |
| | 2018 | 1.325.E+04 | 2.647.E+08 | 1.227.E+00 | 2.750.E-01 | 1.765.E+00 |
| | 2019 | 1.119.E+03 | 7.329.E+06 | 2.419.E+00 | 3.750.E-01 | 4.305.E+00 |
| | 2020 | 1.518.E+03 | 2.996.E+06 | 1.140.E+00 | 2.750.E-01 | 1.756.E+00 |
| | 2021 | 1.117.E+03 | 2.422.E+06 | 1.394.E+00 | 3.250.E-01 | 2.137.E+00 |

|  | Year | | | | | |
|---|---|---|---|---|---|---|
| | 2022 | 3.827.E+03 | 9.639.E+07 | 2.566.E+00 | 3.000.E-01 | 5.616.E+00 |
| | 2023 | 2.217.E+03 | 5.685.E+06 | 1.076.E+00 | 2.917.E-01 | 2.019.E+00 |
| | 2024 | 1.435.E+03 | 3.040.E+06 | 1.215.E+00 | 3.023.E-01 | 1.371.E+00 |
| Tone River | 2012 | 2.475.E+03 | 2.781.E+07 | 2.131.E+00 | 3.171.E-01 | 2.834.E+00 |
| | 2013 | 4.801.E+02 | 2.790.E+05 | 1.100.E+00 | 3.056.E-01 | 1.630.E+00 |
| | 2014 | 2.037.E+03 | 9.273.E+06 | 1.495.E+00 | 3.500.E-01 | 2.769.E+00 |
| | 2015 | 3.099.E+02 | 2.970.E+05 | 1.759.E+00 | 3.429.E-01 | 2.453.E+00 |
| | 2016 | 3.719.E+02 | 3.499.E+05 | 1.591.E+00 | 2.439.E-01 | 2.109.E+00 |
| | 2017 | 8.936.E+02 | 2.024.E+06 | 1.592.E+00 | 3.000.E-01 | 1.930.E+00 |
| | 2018 | 1.728.E+03 | 2.789.E+06 | 9.664.E-01 | 2.927.E-01 | 6.718.E-01 |
| | 2019 | 3.843.E+02 | 2.841.E+05 | 1.387.E+00 | 2.927.E-01 | 2.173.E+00 |
| | 2020 | 1.439.E+03 | 5.706.E+06 | 1.660.E+00 | 3.590.E-01 | 2.459.E+00 |
| | 2021 | 3.897.E+02 | 1.018.E+06 | 2.590.E+00 | 3.333.E-01 | 4.633.E+00 |
| | 2022 | 4.624.E+03 | 1.351.E+08 | 2.514.E+00 | 2.683.E-01 | 4.426.E+00 |
| | 2023 | 2.085.E+03 | 7.952.E+06 | 1.352.E+00 | 2.683.E-01 | 2.345.E+00 |
| | 2024 | 3.732.E+02 | 1.848.E+05 | 1.152.E+00 | 3.171.E-01 | 1.420.E+00 |

Table S3. Identified parameter values with $w=1$ being fixed.

| River | Year | $\alpha$ (-) | $\beta$ (1/day) | $H$ (if $0<H<1$) | $b$ (ind/day) | $\lambda$ (day/ind) | $\mu$ (-) |
|---|---|---|---|---|---|---|---|
| Yahagi River | 1998 | 9.281.E+04 | 6.759.E-06 | | 1.961.E+04 | 9.743.E-06 | 6.024.E-06 |
| | 1999 | 6.469.E+00 | 5.798.E+02 | | 1.110.E+07 | 8.567.E-07 | 2.538.E-08 |
| | 2000 | 5.321.E+01 | 1.481.E+03 | | 1.173.E+08 | 5.872.E-07 | 2.668.E-09 |
| | 2001 | 1.721.E+01 | 1.689.E+03 | | 5.859.E+06 | 5.785.E-06 | 2.041.E-08 |
| | 2002 | 3.534.E+00 | 1.317.E+03 | | 2.749.E+07 | 4.689.E-07 | 9.592.E-09 |
| | 2003 | 6.971.E+04 | 2.951.E-05 | | 2.525.E+03 | 5.296.E-05 | 3.692.E-05 |
| | 2004 | 1.691.E+00 | 3.574.E+01 | 0.655 | 3.421.E+05 | 9.873.E-07 | 4.007.E-07 |
| | 2005 | 3.343.E+04 | 7.586.E-05 | | 7.729.E+03 | 3.507.E-05 | 1.750.E-05 |
| | 2006 | 8.959.E+00 | 2.700.E-01 | | 9.129.E+03 | 1.976.E-05 | 1.222.E-05 |
| | 2007 | 2.699.E+00 | 1.391.E+00 | 0.151 | 6.137.E+04 | 3.047.E-06 | 1.864.E-06 |
| | 2008 | 1.066.E+00 | 1.598.E+15 | 0.967 | 1.005.E+18 | 2.153.E-12 | 1.971.E-19 |
| | 2009 | 5.453.E+05 | 2.648.E-06 | | 1.593.E+04 | 1.255.E-05 | 6.276.E-06 |
| | 2010 | 1.358.E+00 | 2.789.E+02 | 0.821 | 4.894.E+05 | 5.286.E-06 | 5.619.E-07 |
| | 2011 | 2.076.E+00 | 1.852.E+00 | 0.462 | 1.006.E+04 | 3.108.E-05 | 1.596.E-05 |
| | 2012 | 2.025.E+01 | 5.699.E-02 | | 3.854.E+03 | 6.183.E-05 | 3.589.E-05 |
| | 2013 | 9.294.E+04 | 1.575.E-05 | | 5.232.E+03 | 4.818.E-05 | 2.603.E-05 |
| | 2014 | 8.845.E+00 | 1.752.E-01 | | 4.140.E+03 | 5.302.E-05 | 2.991.E-05 |
| | 2015 | 1.009.E+01 | 1.247.E-01 | | 6.787.E+03 | 2.912.E-05 | 1.737.E-05 |
| | 2016 | 1.407.E+00 | 2.538.E+00 | 0.797 | 5.804.E+04 | 3.935.E-06 | 1.930.E-06 |
| | 2017 | 8.130.E+03 | 2.740.E-04 | | 2.491.E+04 | 1.480.E-05 | 5.609.E-06 |
| | 2018 | 3.605.E+00 | 5.140.E-01 | | 1.725.E+04 | 1.305.E-05 | 7.586.E-06 |
| | 2019 | 3.871.E+00 | 1.227.E+03 | | 2.938.E+07 | 9.197.E-07 | 1.059.E-08 |
| | 2020 | 1.188.E+05 | 4.532.E-06 | | 3.930.E+03 | 3.031.E-05 | 1.347.E-05 |
| | 2021 | 1.906.E+00 | 9.587.E+00 | 0.547 | 9.533.E+04 | 7.028.E-06 | 2.045.E-06 |
| | 2022 | 2.335.E+00 | 1.086.E+00 | 0.333 | 7.655.E+03 | 2.701.E-05 | 1.328.E-05 |
| | 2023 | 1.132.E+05 | 9.761.E-06 | | 5.036.E+02 | 3.020.E-04 | 1.679.E-04 |
| Nagara River | 2023 | 1.438.E+00 | 1.053.E+01 | 0.781 | 2.071.E+04 | 2.130.E-05 | 8.190.E-06 |
| | 2024 | 1.592.E+00 | 2.067.E+00 | 0.704 | 6.021.E+03 | 3.034.E-05 | 1.777.E-05 |
| Tama River | 2011 | 2.890.E+00 | 1.656.E+00 | 0.055 | 9.534.E+03 | 2.741.E-05 | 1.437.E-05 |
| | 2012 | 1.068.E+00 | 1.072.E+06 | 0.966 | 7.252.E+08 | 1.540.E-07 | 4.790.E-10 |
| | 2013 | 1.277.E+00 | 1.243.E+01 | 0.862 | 1.165.E+04 | 5.425.E-05 | 1.681.E-05 |
| | 2014 | 1.286.E+00 | 2.947.E+01 | 0.857 | 2.613.E+04 | 4.101.E-05 | 8.814.E-06 |
| | 2015 | 1.082.E+00 | 1.441.E+06 | 0.959 | 3.903.E+08 | 3.857.E-07 | 8.754.E-10 |
| | 2016 | 2.065.E+00 | 1.777.E+00 | 0.468 | 3.596.E+03 | 8.308.E-05 | 3.955.E-05 |
| | 2017 | 1.113.E+00 | 3.695.E+07 | 0.944 | 4.697.E+09 | 1.612.E-07 | 6.441.E-11 |
| | 2018 | 6.010.E+04 | 4.571.E-05 | | 1.110.E+04 | 2.160.E-05 | 1.103.E-05 |
| | 2019 | 1.608.E+04 | 1.401.E-04 | | 2.184.E+03 | 7.370.E-05 | 4.627.E-05 |
| | 2020 | 5.994.E+04 | 2.449.E-05 | | 2.541.E+02 | 8.338.E-04 | 4.716.E-04 |
| | 2021 | 1.797.E+00 | 5.374.E+00 | 0.602 | 7.126.E+02 | 6.823.E-04 | 2.548.E-04 |
| | 2022 | 6.187.E+04 | 1.512.E-05 | | 9.008.E+02 | 2.722.E-04 | 1.306.E-04 |
| | 2023 | 2.217.E+00 | 3.008.E+00 | 0.392 | 3.468.E+03 | 1.584.E-04 | 5.398.E-05 |
| | 2024 | 2.424.E+00 | 2.150.E+00 | 0.288 | 7.305.E+02 | 6.950.E-04 | 2.195.E-04 |
| Ara River | 2012 | 1.537.E+00 | 1.174.E+01 | 0.732 | 8.944.E+04 | 7.579.E-06 | 2.196.E-06 |
| | 2013 | 2.033.E+05 | 5.900.E-06 | | 6.689.E+03 | 4.583.E-05 | 2.267.E-05 |
| | 2014 | 3.787.E+04 | 3.001.E-05 | | 6.022.E+03 | 5.362.E-05 | 2.583.E-05 |
| | 2015 | 1.521.E+05 | 7.020.E-06 | | 1.964.E+03 | 2.326.E-04 | 9.565.E-05 |
| | 2016 | 7.099.E+04 | 2.029.E-05 | | 3.227.E+03 | 1.419.E-04 | 6.148.E-05 |
| | 2017 | 1.138.E+05 | 1.418.E-05 | | 5.611.E+03 | 9.649.E-05 | 3.620.E-05 |
| | 2018 | 1.102.E+04 | 8.251.E-05 | | 6.212.E+03 | 4.708.E-05 | 2.282.E-05 |
| | 2019 | 6.704.E+04 | 3.209.E-05 | | 9.113.E+02 | 2.507.E-04 | 1.558.E-04 |
| | 2020 | 9.420.E+04 | 1.159.E-05 | | 9.421.E+02 | 3.846.E-04 | 1.660.E-04 |
| | 2021 | 6.725.E+00 | 4.529.E+01 | | 2.755.E+05 | 2.332.E-05 | 1.123.E-06 |
| | 2022 | 3.683.E+03 | 1.388.E-03 | | 1.061.E+04 | 3.351.E-05 | 1.534.E-05 |
| | 2023 | 2.491.E+04 | 9.897.E-05 | | 3.782.E+03 | 1.734.E-04 | 5.338.E-05 |
| | 2024 | 2.025.E+03 | 1.165.E-03 | | 2.197.E+03 | 2.549.E-04 | 8.936.E-05 |
| Tone River | 2012 | 7.366.E+01 | 2.226.E+01 | | 3.886.E+06 | 2.694.E-06 | 7.919.E-08 |
| | 2013 | 1.729.E+00 | 1.189.E+01 | 0.636 | 3.383.E+03 | 3.942.E-04 | 7.348.E-05 |
| | 2014 | 2.401.E+01 | 2.991.E+00 | | 1.261.E+05 | 2.470.E-05 | 2.496.E-06 |
| | 2015 | 4.707.E+01 | 6.619.E+00 | | 8.907.E+04 | 6.337.E-05 | 3.629.E-06 |
| | 2016 | 1.148.E+01 | 1.289.E+02 | | 4.915.E+05 | 2.297.E-05 | 4.857.E-07 |
| | 2017 | 9.473.E+04 | 1.449.E-05 | | 5.923.E+02 | 4.729.E-04 | 2.446.E-04 |
| | 2018 | 7.162.E+01 | 2.004.E+01 | | 2.412.E+06 | 8.671.E-06 | 1.197.E-07 |
| | 2019 | 5.787.E+01 | 1.137.E+01 | | 2.414.E+05 | 4.050.E-05 | 1.177.E-06 |
| | 2020 | 3.544.E+04 | 8.411.E-05 | | 2.430.E+03 | 1.930.E-04 | 8.369.E-05 |
| | 2021 | 1.102.E+05 | 6.060.E-04 | | 2.169.E+04 | 7.668.E-05 | 1.280.E-05 |
| | 2022 | 4.934.E+01 | 9.344.E+00 | | 1.964.E+06 | 2.162.E-06 | 1.285.E-07 |
| | 2023 | 4.743.E+04 | 3.429.E-05 | | 1.970.E+03 | 1.890.E-04 | 7.914.E-05 |
| | 2024 | 5.590.E+01 | 9.987.E+00 | | 1.990.E+05 | 5.672.E-05 | 1.549.E-06 |

Table S4. Identified parameter values with $w$ being fitted (magenta indicates the Negative $b$ values).

| River | Year | $\alpha$ (-) | $\beta$ (1/day) | $b$ (ind/day) | $\lambda$ (day/ind) | $\mu$ (-) | $w$ (-) |
|---|---|---|---|---|---|---|---|
| Yahagi River | 1998 | 9.281.E+04 | 2.579.E-06 | -1.216.E+04 | 9.743.E-06 | 1.579.E-05 | 0.000.E+00 |
| | 1999 | 6.469.E+00 | 5.626.E+02 | 1.076.E+07 | 8.567.E-07 | 2.616.E-08 | 0.000.E+00 |
| | 2000 | 5.321.E+01 | 1.475.E+03 | 1.167.E+08 | 5.872.E-07 | 2.680.E-09 | 0.000.E+00 |
| | 2001 | 1.721.E+01 | 1.683.E+03 | 5.838.E+06 | 5.785.E-06 | 2.048.E-08 | 0.000.E+00 |
| | 2002 | 3.534.E+00 | 1.290.E+03 | 2.692.E+07 | 4.689.E-07 | 9.793.E-09 | 0.000.E+00 |
| | 2003 | 6.971.E+04 | 2.387.E-05 | 9.323.E+02 | 5.296.E-05 | 4.564.E-05 | 7.258.E-01 |
| | 2004 | 1.691.E+00 | 2.123.E+01 | 1.084.E+05 | 9.873.E-07 | 6.745.E-07 | 0.000.E+00 |
| | 2005 | 3.343.E+04 | 3.800.E-05 | 3.041.E+01 | 3.507.E-05 | 3.493.E-05 | 0.000.E+00 |
| | 2006 | 8.959.E+00 | 1.543.E-01 | -1.123.E+03 | 1.976.E-05 | 2.139.E-05 | 3.069.E-01 |
| | 2007 | 2.699.E+00 | 1.079.E+00 | 2.586.E+04 | 3.047.E-06 | 2.404.E-06 | 6.330.E-01 |
| | 2008 | 1.066.E+00 | 1.598.E+15 | 1.005.E+18 | 2.153.E-12 | 1.971.E-19 | 1.000.E+00 |
| | 2009 | 5.453.E+05 | 1.323.E-06 | -8.292.E+00 | 1.255.E-05 | 1.256.E-05 | 0.000.E+00 |
| | 2010 | 1.358.E+00 | 2.492.E+02 | 4.311.E+05 | 5.286.E-06 | 6.288.E-07 | 0.000.E+00 |
| | 2011 | 2.076.E+00 | 9.009.E-01 | -5.566.E+02 | 3.108.E-05 | 3.280.E-05 | 0.000.E+00 |
| | 2012 | 2.025.E+01 | 2.391.E-02 | -1.478.E+03 | 6.183.E-05 | 8.553.E-05 | 0.000.E+00 |
| | 2013 | 9.294.E+04 | 1.547.E-05 | 5.029.E+03 | 4.818.E-05 | 2.650.E-05 | 9.670.E-01 |
| | 2014 | 8.845.E+00 | 9.394.E-02 | -2.657.E+02 | 5.302.E-05 | 5.579.E-05 | 1.778.E-01 |
| | 2015 | 1.009.E+01 | 5.856.E-02 | -2.132.E+03 | 2.912.E-05 | 3.698.E-05 | 1.109.E-01 |
| | 2016 | 1.407.E+00 | 1.293.E+00 | 2.186.E+03 | 3.935.E-06 | 3.787.E-06 | 0.000.E+00 |
| | 2017 | 8.130.E+03 | 1.702.E-04 | 9.707.E+03 | 1.480.E-05 | 9.032.E-06 | 0.000.E+00 |
| | 2018 | 3.605.E+00 | 2.426.E-01 | -4.503.E+03 | 1.305.E-05 | 1.607.E-05 | 9.172.E-02 |
| | 2019 | 3.871.E+00 | 1.213.E+03 | 2.904.E+07 | 9.197.E-07 | 1.071.E-08 | 0.000.E+00 |
| | 2020 | 1.188.E+05 | 2.518.E-06 | 7.879.E+02 | 3.031.E-05 | 2.423.E-05 | 0.000.E+00 |
| | 2021 | 1.906.E+00 | 6.797.E+00 | 5.620.E+04 | 7.028.E-06 | 2.885.E-06 | 0.000.E+00 |
| | 2022 | 2.335.E+00 | 5.521.E-01 | 2.506.E+02 | 2.701.E-05 | 2.613.E-05 | 0.000.E+00 |
| | 2023 | 1.132.E+05 | 4.334.E-06 | -1.270.E+02 | 3.020.E-04 | 3.782.E-04 | 0.000.E+00 |
| Nagara River | 2023 | 1.438.E+00 | 6.484.E+00 | 7.775.E+03 | 2.130.E-05 | 1.331.E-05 | 0.000.E+00 |
| | 2024 | 1.592.E+00 | 8.564.E-01 | -2.521.E+03 | 3.045.E-05 | 4.332.E-05 | 0.000.E+00 |
| Tama River | 2011 | 2.890.E+00 | 1.129.E+00 | 3.164.E+03 | 2.741.E-05 | 2.107.E-05 | 3.935.E-01 |
| | 2012 | 1.068.E+00 | 1.069.E+06 | 7.229.E+08 | 1.540.E-07 | 4.805.E-10 | 0.000.E+00 |
| | 2013 | 1.277.E+00 | 8.578.E+00 | 6.417.E+03 | 5.425.E-05 | 2.436.E-05 | 0.000.E+00 |
| | 2014 | 1.286.E+00 | 2.314.E+01 | 1.897.E+04 | 4.101.E-05 | 1.123.E-05 | 0.000.E+00 |
| | 2015 | 1.082.E+00 | 1.438.E+06 | 3.895.E+08 | 3.857.E-07 | 8.773.E-10 | 0.000.E+00 |
| | 2016 | 2.065.E+00 | 9.311.E-01 | 3.294.E+02 | 8.308.E-05 | 7.548.E-05 | 0.000.E+00 |
| | 2017 | 1.113.E+00 | 3.693.E+07 | 4.695.E+09 | 1.612.E-07 | 6.443.E-11 | 0.000.E+00 |
| | 2018 | 6.010.E+04 | 3.651.E-05 | 6.540.E+03 | 2.160.E-05 | 1.380.E-05 | 6.059.E-01 |
| | 2019 | 1.608.E+04 | 1.142.E-04 | 1.098.E+03 | 7.370.E-05 | 5.677.E-05 | 7.052.E-01 |
| | 2020 | 5.994.E+04 | 1.260.E-05 | -2.995.E+01 | 8.338.E-04 | 9.168.E-04 | 1.415.E-01 |
| | 2021 | 1.797.E+00 | 3.367.E+00 | 2.877.E+02 | 6.823.E-04 | 4.068.E-04 | 0.000.E+00 |
| | 2022 | 6.187.E+04 | 7.866.E-06 | 6.979.E+01 | 2.722.E-04 | 2.511.E-04 | 0.000.E+00 |
| | 2023 | 2.217.E+00 | 1.983.E+00 | 1.676.E+03 | 1.584.E-04 | 8.187.E-05 | 0.000.E+00 |
| | 2024 | 2.424.E+00 | 1.471.E+00 | 3.933.E+02 | 6.950.E-04 | 3.208.E-04 | 0.000.E+00 |
| Ara River | 2012 | 1.537.E+00 | 8.334.E+00 | 5.294.E+04 | 7.579.E-06 | 3.093.E-06 | 0.000.E+00 |
| | 2013 | 2.033.E+05 | 2.982E-06 | 1.433.E+02 | 4.583.E-05 | 4.485.E-05 | 0.000.E+00 |
| | 2014 | 3.787.E+04 | 1.556E-05 | 4.269.E+02 | 5.362.E-05 | 4.982.E-05 | 0.000.E+00 |
| | 2015 | 1.521.E+05 | 4.133E-06 | 5.921.E+02 | 2.326.E-04 | 1.625.E-04 | 0.000.E+00 |
| | 2016 | 7.099.E+04 | 1.150E-05 | 7.602.E+02 | 1.419.E-04 | 1.085.E-04 | 0.000.E+00 |
| | 2017 | 1.138.E+05 | 8.859E-06 | 2.242.E+03 | 9.649.E-05 | 5.793.E-05 | 0.000.E+00 |
| | 2018 | 1.102.E+04 | 4.253E-05 | 3.710.E+02 | 4.708.E-05 | 4.427.E-05 | 0.000.E+00 |
| | 2019 | 6.704.E+04 | 2.253E-05 | 1.940.E+02 | 2.507.E-04 | 2.219.E-04 | 5.205.E-01 |
| | 2020 | 9.420.E+04 | 6.587E-06 | 2.271.E+02 | 3.846.E-04 | 2.919.E-04 | 0.000.E+00 |
| | 2021 | 6.725.E+00 | 4.311E+01 | 2.616.E+05 | 2.332.E-05 | 1.179.E-06 | 0.000.E+00 |
| | 2022 | 3.683.E+03 | 7.526E-04 | 1.653.E+03 | 3.351.E-05 | 2.829.E-05 | 0.000.E+00 |
| | 2023 | 2.491.E+04 | 6.851E-05 | 2.100.E+03 | 1.734.E-04 | 7.711.E-05 | 0.000.E+00 |

|            | Year | | | | | | |
|------------|------|-----------|-----------|-----------|-----------|-----------|-----------|
|            | 2024 | 2.025.E+03 | 7.563E-04 | 1.011.E+03 | 2.549.E-04 | 1.376.E-04 | 0.000.E+00 |
|            | 2012 | 7.366.E+01 | 2.161.E+01 | 3.769.E+06 | 2.694.E-06 | 8.159.E-08 | 0.000.E+00 |
|            | 2013 | 1.729.E+00 | 9.670.E+00 | 2.608.E+03 | 3.942.E-04 | 9.032.E-05 | 0.000.E+00 |
|            | 2014 | 2.401.E+01 | 2.689.E+00 | 1.119.E+05 | 2.470.E-05 | 2.776.E-06 | 0.000.E+00 |
|            | 2015 | 4.707.E+01 | 6.240.E+00 | 8.366.E+04 | 6.337.E-05 | 3.849.E-06 | 0.000.E+00 |
|            | 2016 | 1.148.E+01 | 1.262.E+02 | 4.809.E+05 | 2.297.E-05 | 4.962.E-07 | 0.000.E+00 |
|            | 2017 | 9.473.E+04 | 6.997.E-06 | -4.213.E+01 | 4.729.E-04 | 5.065.E-04 | 0.000.E+00 |
| Tone River | 2018 | 7.162.E+01 | 1.976.E+01 | 2.378.E+06 | 8.671.E-06 | 1.214.E-07 | 0.000.E+00 |
|            | 2019 | 5.787.E+01 | 1.104.E+01 | 2.342.E+05 | 4.050.E-05 | 1.213.E-06 | 0.000.E+00 |
|            | 2020 | 3.544.E+04 | 4.764.E-05 | 5.700.E+02 | 1.930.E-04 | 1.477.E-04 | 0.000.E+00 |
|            | 2021 | 1.103.E+05 | 5.048.E-04 | 1.735.E+04 | 7.668.E-05 | 1.537.E-05 | 0.000.E+00 |
|            | 2022 | 4.934.E+01 | 8.789.E+00 | 1.840.E+06 | 2.162.E-06 | 1.366.E-07 | 0.000.E+00 |
|            | 2023 | 4.743.E+04 | 1.993.E-05 | 5.505.E+02 | 1.890.E-04 | 1.362.E-04 | 0.000.E+00 |
|            | 2024 | 5.590.E+01 | 9.714.E+00 | 1.935.E+05 | 5.672.E-05 | 1.593.E-06 | 0.000.E+00 |

**Table S5.** Empirical skewness, with $w$ being fixed and fitted. Their relative errors are presented as well. "Improvement" in the extreme right column is the relative error of $w$ being fitted divided by that of $w$ being fixed as 1; a small value implies a substantial improvement in the reproduction of skewness. Magenta indicates cases with negative $b$ values in **Table S4**.

| River | Year | $w$ is fixed to 1 | | $w$ is fitted | | Improvement |
|---|---|---|---|---|---|---|
| | | Skw (-) | Relative error | Skw (-) | Relative error | |
| Yahagi River | 1998 | 2.846.E+00 | 8.367.E-01 | 1.759.E+00 | 1.349.E-01 | 1.613.E-01 |
| | 1999 | 2.120.E+02 | 4.180.E+01 | 2.059.E+02 | 4.056.E+01 | 9.705.E-01 |
| | 2000 | 9.943.E+02 | 3.386.E+02 | 9.898.E+02 | 3.371.E+02 | 9.955.E-01 |
| | 2001 | 9.568.E+02 | 5.660.E+02 | 9.535.E+02 | 5.640.E+02 | 9.965.E-01 |
| | 2002 | 2.249.E+02 | 6.612.E+01 | 2.204.E+02 | 6.478.E+01 | 9.797.E-01 |
| | 2003 | 4.831.E+00 | 4.404.E-02 | 4.627.E+00 | 9.876.E-10 | 2.242.E-08 |
| | 2004 | 2.242.E+01 | 2.478.E+00 | 1.595.E+01 | 1.474.E+00 | 5.948.E-01 |
| | 2005 | 6.504.E+00 | 1.272.E+00 | 4.339.E+00 | 5.155.E-01 | 4.053.E-01 |
| | 2006 | 5.421.E+00 | 2.149.E-01 | 4.462.E+00 | 2.005.E-07 | 9.331.E-07 |
| | 2007 | 5.688.E+00 | 7.497.E-02 | 5.292.E+00 | 1.148.E-07 | 1.531.E-06 |
| | 2008 | 4.620.E+07 | 1.375.E+07 | 4.620.E+07 | 1.375.E+07 | 1.000.E+00 |
| | 2009 | 5.687.E+00 | 1.242.E+00 | 3.799.E+00 | 4.981.E-01 | 4.009.E-01 |
| | 2010 | 3.768.E+01 | 1.025.E+01 | 3.406.E+01 | 9.169.E+00 | 8.945.E-01 |
| | 2011 | 5.189.E+00 | 5.787.E-01 | 3.429.E+00 | 4.313.E-02 | 7.452.E-02 |
| | 2012 | 3.735.E+00 | 6.209.E-01 | 2.363.E+00 | 2.562.E-02 | 4.126.E-02 |
| | 2013 | 4.643.E+00 | 5.445.E-03 | 4.618.E+00 | 1.625.E-08 | 2.985.E-06 |
| | 2014 | 4.543.E+00 | 3.177.E-01 | 3.448.E+00 | 9.551.E-08 | 3.006.E-07 |
| | 2015 | 3.995.E+00 | 3.993.E-01 | 2.855.E+00 | 5.062.E-08 | 1.268.E-07 |
| | 2016 | 4.614.E+00 | 1.945.E+00 | 3.096.E+00 | 9.757.E-01 | 5.018.E-01 |
| | 2017 | 6.839.E+00 | 1.682.E+00 | 4.959.E+00 | 9.450.E-01 | 5.618.E-01 |
| | 2018 | 4.236.E+00 | 4.208.E-01 | 2.981.E+00 | 1.005.E-07 | 2.389.E-07 |
| | 2019 | 2.128.E+02 | 8.999.E+01 | 2.104.E+02 | 8.895.E+01 | 9.885.E-01 |
| | 2020 | 5.119.E+00 | 1.762.E+00 | 3.545.E+00 | 9.124.E-01 | 5.178.E-01 |
| | 2021 | 1.222.E+01 | 3.174.E+00 | 9.465.E+00 | 2.233.E+00 | 7.036.E-01 |
| | 2022 | 5.727.E+00 | 9.743.E-01 | 3.840.E+00 | 3.236.E-01 | 3.321.E-01 |
| | 2023 | 4.994.E+00 | 7.891.E-01 | 3.209.E+00 | 1.498.E-01 | 1.899.E-01 |
| Nagara River | 2023 | 4.994.E+00 | 7.891.E-01 | 3.209.E+00 | 1.498.E-01 | 1.899.E-01 |
| | 2024 | 4.925.E+00 | 1.990.E+00 | 3.439.E+00 | 1.088.E+00 | 5.467.E-01 |
| Tama River | 2011 | 8.892.E+00 | 2.328.E+00 | 6.423.E+00 | 1.404.E+00 | 6.030.E-01 |
| | 2012 | 4.441.E+00 | 6.467.E-01 | 2.801.E+00 | 3.215.E-02 | 4.972.E-02 |
| | 2013 | 6.933.E+00 | 1.703.E-01 | 5.924.E+00 | 5.087.E-08 | 2.987.E-07 |
| | 2014 | 9.189.E+02 | 4.095.E+02 | 9.160.E+02 | 4.082.E+02 | 9.969.E-01 |
| | 2015 | 7.581.E+00 | 3.274.E+00 | 5.787.E+00 | 2.263.E+00 | 6.911.E-01 |
| | 2016 | 1.154.E+01 | 4.613.E+00 | 9.501.E+00 | 3.620.E+00 | 7.848.E-01 |
| | 2017 | 1.177.E+03 | 7.774.E+02 | 1.174.E+03 | 7.756.E+02 | 9.977.E-01 |
| | 2018 | 5.643.E+00 | 9.531.E-01 | 3.823.E+00 | 3.232.E-01 | 3.391.E-01 |
| | 2019 | 7.432.E+03 | 2.082.E+03 | 7.429.E+03 | 2.081.E+03 | 9.996.E-01 |
| | 2020 | 7.006.E+00 | 8.886.E-02 | 6.434.E+00 | 1.262.E-12 | 1.420.E-11 |
| | 2021 | 5.722.E+00 | 5.477.E-02 | 5.425.E+00 | 1.593.E-08 | 2.909.E-07 |
| | 2022 | 4.761.E+00 | 3.549.E-01 | 3.514.E+00 | 2.391.E-07 | 6.737.E-07 |
| | 2023 | 8.360.E+00 | 1.372.E+00 | 6.087.E+00 | 7.272.E-01 | 5.299.E-01 |
| | 2024 | 4.341.E+00 | 1.148.E+00 | 2.934.E+00 | 4.517.E-01 | 3.934.E-01 |
| Ara River | 2012 | 1.038.E+01 | 2.196.E+00 | 8.046.E+00 | 1.478.E+00 | 6.730.E-01 |
| | 2013 | 4.250.E+00 | 1.188.E+00 | 2.843.E+00 | 4.638.E-01 | 3.905.E-01 |
| | 2014 | 4.153.E+00 | 1.106.E+00 | 2.803.E+00 | 4.215.E-01 | 3.811.E-01 |
| | 2015 | 3.962.E+00 | 2.462.E+00 | 2.807.E+00 | 1.453.E+00 | 5.902.E-01 |
| | 2016 | 4.378.E+00 | 1.291.E+00 | 3.054.E+00 | 5.986.E-01 | 4.636.E-01 |
| | 2017 | 4.845.E+00 | 1.732.E+00 | 3.523.E+00 | 9.867.E-01 | 5.697.E-01 |

|  | Year | | | | | |
|---|---|---|---|---|---|---|
| | 2018 | 3.876.E+00 | 1.196.E+00 | 2.611.E+00 | 4.795.E-01 | 4.007.E-01 |
| | 2019 | 4.779.E+00 | 1.100.E-01 | 4.305.E+00 | 2.743.E-08 | 2.494.E-07 |
| | 2020 | 4.300.E+00 | 1.449.E+00 | 3.004.E+00 | 7.106.E-01 | 4.905.E-01 |
| | 2021 | 5.777.E+01 | 2.603.E+01 | 5.511.E+01 | 2.478.E+01 | 9.523.E-01 |
| | 2022 | 8.861.E+00 | 5.780.E-01 | 6.079.E+00 | 8.250.E-02 | 1.427.E-01 |
| | 2023 | 6.326.E+00 | 2.134.E+00 | 4.837.E+00 | 1.396.E+00 | 6.543.E-01 |
| | 2024 | 6.078.E+00 | 3.433.E+00 | 4.500.E+00 | 2.282.E+00 | 6.648.E-01 |
| Tone River | 2012 | 1.449.E+02 | 5.013.E+01 | 1.407.E+02 | 4.867.E+01 | 9.709.E-01 |
| | 2013 | 1.140.E+01 | 5.990.E+00 | 9.605.E+00 | 4.892.E+00 | 8.167.E-01 |
| | 2014 | 2.928.E+01 | 9.576.E+00 | 2.659.E+01 | 8.605.E+00 | 8.986.E-01 |
| | 2015 | 6.123.E+01 | 2.396.E+01 | 5.791.E+01 | 2.261.E+01 | 9.436.E-01 |
| | 2016 | 1.503.E+02 | 7.029.E+01 | 1.472.E+02 | 6.882.E+01 | 9.790.E-01 |
| | 2017 | 4.510.E+00 | 1.337.E+00 | 2.973.E+00 | 5.403.E-01 | 4.041.E-01 |
| | 2018 | 1.400.E+02 | 2.074.E+02 | 1.381.E+02 | 2.046.E+02 | 9.863.E-01 |
| | 2019 | 9.534.E+01 | 4.288.E+01 | 9.265.E+01 | 4.164.E+01 | 9.711.E-01 |
| | 2020 | 6.218.E+00 | 1.529.E+00 | 4.337.E+00 | 7.641.E-01 | 4.998.E-01 |
| | 2021 | 3.016.E+01 | 5.510.E+00 | 2.585.E+01 | 4.579.E+00 | 8.310.E-01 |
| | 2022 | 8.432.E+01 | 1.805.E+01 | 7.959.E+01 | 1.698.E+01 | 9.408.E-01 |
| | 2023 | 5.326.E+00 | 1.271.E+00 | 3.754.E+00 | 6.007.E-01 | 4.726.E-01 |
| | 2024 | 8.426.E+01 | 5.833.E+01 | 8.202.E+01 | 5.675.E+01 | 9.730.E-01 |

**Table S6.** $R^2$ values between each pair among the Yahagi, Tama, Ara, and Tone Rivers between 2012 and 2023.

|  | Yahagi River | Tama River | Ara River | Tone River |
|---|---|---|---|---|
| Yahagi River | 1 | 0.0082 | 0.0066 | 0.0679 |
| Tama River |  | 1 | 0.9272 | 0.4268 |
| Ara River |  |  | 1 | 0.0862 |
| Tone River |  |  |  | 1 |

*Correction note to: Yoshioka, H. (2025). Chaos, Solitons & Fractals. Vol. 192, 115911.*
*or its preprint version: arXiv:2411.12272v2*

Hidekazu Yoshioka[1,*]

[1] Japan Advanced Institute of Science and Technology, 1-1 Asahidai, Nomi 923-1292, Japan
[*] Corresponding author: yoshih@jaist.ac.jp, ORCID: 0000-0002-5293-3246

**Abstract**
This is a correction note to Proposition 3 in the following paper (called Y25 in this letter): **Yoshioka H. (2025). Superposition of interacting stochastic processes with memory and its application to migrating fish counts. Chaos, Solitons & Fractals. Vol. 192, 115911. https://doi.org/10.1016/j.chaos.2024.115911 (or its preprint version https://arxiv.org/abs/2411.12272)**

**Main text**
The paper Y25[1] (Chaos, Solitons & Fractals, Vol. 192, 115911) discussed a superposition process model, a non-Markovian stochastic process model, with its application to fish migration. A generalized Riccati equation associated with the model was also studied in the paper as a separate topic (Proposition 3 [1]). The discussion and correction below do not concern the superposition model itself and its application, but only the generalized Riccati equation.

While investigating a generalization of the model mentioned above, we found that Proof of Proposition 3 of Section A.4 [1] was incorrect at least in the following two points.

✓ First, Eq. (63) in Y25[1] was incorrect because it neglects the term proportional to the constant $w$ (e.g., the first integrand of Eq. (60) in Y25[1]).

✓ Second, even if this neglected term is recovered, it is difficult (at least for the author) to obtain a Lipschitz continuity as concluded in the last line of this equation because of the nonlinearity of the neglected term and proportionality with respect to the unbounded reversion speed $r$. More specifically, the corresponding term seems not to be controllable in the Lebesgue space of integrable functions with respect to probability measures considered in Y25[1].

Note that these issues are not encountered for the Lyapunov equations (Proposition 2 [1]) due to that they are (fortunately) linear and do not have corresponding terms. In the sequel, we focus on the second issue because the first one is resolved once the second one is.

In this letter, we present a corrected version of Proposition 3 [1] and its proof. We show that their Riccati equation can be well-defined in a different setting, where solutions to the equation are defined in a space of bounded continuous functions instead of that of integrable functions that was originally considered in Y25[1]. The streamline of the new proof is similar to that of the previous one but uses different functional spaces to which solutions belong.

By a generalized Riccati equation, we mean the following initial value problem to find a function $\varphi:[0,+\infty)\times(0,+\infty)\to\mathbb{R}$ (understood as an ordinary differential equation in a Banach space):

$$\frac{\mathrm{d}\varphi(t,r)}{\mathrm{d}t} = -r\varphi(t,r) + wr\int_{z=0}^{z=+\infty}\left(1-e^{-\varphi(t,r)z}\right)v(\mathrm{d}z) \\ + (1-w)r\int_{u=0}^{u=+\infty}\int_{z=0}^{z=+\infty}\left(1-e^{-\varphi(t,u)z}\right)v(\mathrm{d}z)\rho(\mathrm{d}u), \quad t>0,\ r>0 \quad (83)$$

subject to the initial condition $\varphi(0,\cdot) = \theta > 0$, and $w\in[0,1]$ is a constant. A mild form of (83) is set as

$$\varphi(t,r) = \theta e^{-rt} + \int_{s=0}^{s=t} e^{-r(t-s)} r \left\{ w\int_{z=0}^{z=+\infty}\left(1-e^{-\varphi(s,r)z}\right)v(\mathrm{d}z) \\ + (1-w)\int_{u=0}^{u=+\infty}\int_{z=0}^{z=+\infty}\left(1-e^{-\varphi(s,u)z}\right)v(\mathrm{d}z)\rho(\mathrm{d}u) \right\} \mathrm{d}s, \quad t\geq 0,\ r>0. \quad (84)$$

$$(:= \mathbb{G}(\varphi)(t,r))$$

Here, $\rho$ is a probability measure of a positive random variable and $\nu$ is a Lévy measure of a pure-jump subordinator such that $M_k = \int_0^{+\infty} z^k \nu(\mathrm{d}z)$ with $M_1 \in (0,1)$. The generalized Riccati equation of the form (83) was presented in Y25[1] as a governing equation of a coefficient in the moment generating function of a superposition process. We assume that $\rho(\mathrm{d}r)$ admits a positive density that is absolutely continuous with respect to the Lebesgue measure $\mathrm{d}r$ (for cases where $\rho$ is singular, see **Remark 3** at the bottom of this letter).

By a solution to the generalized Riccati equation (83), we mean a bounded continuous function in $[0,+\infty) \times (0,+\infty)$ that satisfies (84). More specifically, we set the space of bounded continuous functions equipped with the supremum norm $\|\cdot\|$, which is a Banach space (e.g., p.30 in Clason [2]):

$$C_{\mathrm{b}} = \left\{ \varphi \in C(0,+\infty) : \|\varphi\| := \sup_{r \in (0,+\infty)} |\varphi(r)| < +\infty \right\}. \tag{85}$$

For each fixed $T > 0$, we set another Banach space $C_{\mathrm{b},T}$:

$$C_{\mathrm{b},T} = \left\{ \varphi \in C([0,T] \times (0,+\infty)) : \|\varphi\|_T := \sup_{0 \le t \le T} \|\varphi(t,\cdot)\| < +\infty \right\}. \tag{86}$$

We set a constant $\bar{\varphi}(\theta) = \theta(1 - M_1)^{-1} > 0$ and a truncation function $\hat{\varphi}(\cdot,\cdot) = \max\{0, \min\{\bar{\varphi}(\theta), \varphi(\cdot,\cdot)\}\}$. This $\bar{\varphi}(\theta)$ turns out to be an upper bound of solutions.

With these preparations, the corrected version of Proposition 3 [1] is presented as follows with notations and equation numberings used in that paper.

*A corrected version of Proposition 3*

*For any $T>0$, the generalized Riccati equation (22) subject to an initial condition $B_0(\cdot,\theta)=\theta$ ($\theta \geq 0$) admits a unique nonnegative solution $B \in C_{b,T}$ with $\|B\|_T \leq \theta(1-M_1)^{-1}$.*

We present a proof of this proposition. We first consider an auxiliary equation with the truncation:
$$\varphi(t,r) = \hat{\mathbb{G}}(\varphi)(t,r) = \mathbb{G}(\hat{\varphi})(t,r), \ t \geq 0, \ r > 0. \tag{87}$$

For each $T > 0$, $\hat{\mathbb{G}}$ is a mapping from $C_{b,T}$ to $C_{b,T}$ because $\theta e^{-rt}$ is bounded continuous and the following integrals for each $t \in [0,T]$ and $r > 0$ are bounded continuous as well: $r\int_{s=0}^{s=t} e^{-r(t-s)}\mathrm{d}s = 1 - e^{-rt}$, $r\int_{s=0}^{s=t}e^{-r(t-s)}\int_{z=0}^{z=+\infty}\left(1-e^{-\hat{\varphi}(s,r)z}\right)v(\mathrm{d}z)\mathrm{d}s$, and $r\int_{s=0}^{s=t}e^{-r(t-s)}\int_{u=0}^{u=+\infty}\int_{z=0}^{z=+\infty}\left(1-e^{-\hat{\varphi}(s,u)z}\right)v(\mathrm{d}z)\rho(\mathrm{d}u)\mathrm{d}s$. The continuity follows from dominated convergence theorem (Theorem 4.20 in Farenick [3]) due to that integrands of these integrals are compositions of continuous functions that can be bounded by integrable functions, i.e., $0 \leq e^{-r(t-s)} \leq 1$ and $0 \leq 1 - e^{-\hat{\varphi}(s,u)z} \leq \hat{\varphi}(s,u)z \leq \bar{\varphi}(\theta)z$. The range of integration $(0,t)$ in these integrals can be transformed to $(0,1)$ by introducing the new auxiliary parameter $s' = ts$ without losing the boundedness and continuity.

Fix $T > 0$ and $\varphi_1, \varphi_2 \in C_{b,T}$. We have the boundedness of $\hat{\mathbb{G}}$ as follows:

$$\left\|\hat{\mathbb{G}}(\varphi_1)\right\|_T = \left\|\mathbb{G}(\hat{\varphi}_1)\right\|_T$$
$$= \left\|\theta e^{-rt} + \int_{s=0}^{s=t} e^{-r(t-s)}r\left\{\begin{array}{l}w\int_{z=0}^{z=+\infty}\left(1-e^{-\hat{\varphi}_1(s,r)z}\right)v(\mathrm{d}z)\\+(1-w)\int_{u=0}^{u=+\infty}\int_{z=0}^{z=+\infty}\left(1-e^{-\hat{\varphi}_1(s,u)z}\right)v(\mathrm{d}z)\rho(\mathrm{d}u)\end{array}\right\}\mathrm{d}s\right\|_T$$
$$\leq \theta + w\left\|\int_{s=0}^{s=t}e^{-r(t-s)}r\int_{z=0}^{z=+\infty}\left(1-e^{-\hat{\varphi}_1(s,r)z}\right)v(\mathrm{d}z)\mathrm{d}s\right\|_T$$
$$+(1-w)\left\|\int_{s=0}^{s=t}e^{-r(t-s)}r\int_{u=0}^{u=+\infty}\int_{z=0}^{z=+\infty}\left(1-e^{-\hat{\varphi}_1(s,u)z}\right)v(\mathrm{d}z)\rho(\mathrm{d}u)\mathrm{d}s\right\|_T$$
$$\leq \theta + w\left\|\int_{s=0}^{s=t}e^{-r(t-s)}r\int_{z=0}^{z=+\infty}\left(\hat{\varphi}_1(s,r)z\right)v(\mathrm{d}z)\mathrm{d}s\right\|_T \tag{88}$$
$$+(1-w)\left\|\int_{s=0}^{s=t}e^{-r(t-s)}r\int_{u=0}^{u=+\infty}\int_{z=0}^{z=+\infty}\left(\hat{\varphi}_1(s,u)z\right)v(\mathrm{d}z)\rho(\mathrm{d}u)\mathrm{d}s\right\|_T$$
$$\leq \theta + wM_1\bar{\varphi}(\theta)\left\|\int_{s=0}^{s=t}e^{-r(t-s)}r\mathrm{d}s\right\|_T + (1-w)M_1\bar{\varphi}(q)\left\|\int_{s=0}^{s=t}e^{-r(t-s)}r\int_{u=0}^{u=+\infty}\rho(\mathrm{d}u)\mathrm{d}s\right\|_T$$
$$\leq \theta + wM_1\bar{\varphi}(\theta)\left\|\int_{s=0}^{s=t}e^{-r(t-s)}r\mathrm{d}s\right\|_T + (1-w)M_1\bar{\varphi}(q)\left\|\int_{s=0}^{s=t}e^{-r(t-s)}r\mathrm{d}s\right\|_T$$
$$\leq \theta + M_1\bar{\varphi}(\theta)$$

Here, we used
$$\left\|\int_{s=0}^{s=t}e^{-r(t-s)}r\mathrm{d}s\right\|_T = \sup_{0 \leq t \leq T}\sup_{r>0}\int_{s=0}^{s=t}e^{-r(t-s)}r\mathrm{d}s \leq \sup_{0 \leq t \leq T}\sup_{r>0}\left(1-e^{-rt}\right) \leq 1, \tag{89}$$

which will be used multiple times in the sequel. For the Lipschitz continuity of $\hat{\mathbb{G}}$, we have

$$\left\|\hat{\mathbb{G}}(\varphi_1)-\hat{\mathbb{G}}(\varphi_2)\right\|_T = \left\|\mathbb{G}(\hat{\varphi}_1)-\mathbb{G}(\hat{\varphi}_2)\right\|_T$$

$$= \left\|\int_{s=0}^{s=t} e^{-r(t-s)} r \left\{ w\int_{z=0}^{z=+\infty}\left(e^{-\hat{\varphi}_2(s,r)z}-e^{-\hat{\varphi}_1(s,r)z}\right)v(\mathrm{d}z) \right. \right.$$
$$\left.\left. +(1-w)\int_{u=0}^{u=+\infty}\int_{z=0}^{z=+\infty}\left(e^{-\hat{\varphi}_2(s,u)z}-e^{-\hat{\varphi}_1(s,u)z}\right)v(\mathrm{d}z)\rho(\mathrm{d}u)\right\}\mathrm{d}s\right\|_T. \quad (90)$$

$$\leq w\left\|\int_{s=0}^{s=t}e^{-r(t-s)}r\int_{z=0}^{z=+\infty}\left|e^{-\hat{\varphi}_2(s,r)z}-e^{-\hat{\varphi}_1(s,r)z}\right|v(\mathrm{d}z)\mathrm{d}s\right\|_T$$
$$+(1-w)\left\|\int_{s=0}^{s=t}e^{-r(t-s)}r\int_{u=0}^{u=+\infty}\int_{z=0}^{z=+\infty}\left|e^{-\hat{\varphi}_2(s,u)z}-e^{-\hat{\varphi}_1(s,u)z}\right|v(\mathrm{d}z)\rho(\mathrm{d}u)\mathrm{d}s\right\|_T$$

We also have

$$\left\|\int_{s=0}^{s=t}e^{-r(t-s)}r\int_{z=0}^{z=+\infty}\left|e^{-\hat{\varphi}_2(s,r)z}-e^{-\hat{\varphi}_1(s,r)z}\right|v(\mathrm{d}z)\mathrm{d}s\right\|_T$$
$$\leq M_1\left\|\int_{s=0}^{s=t}e^{-r(t-s)}r\left|\hat{\varphi}_2(s,r)-\hat{\varphi}_1(s,r)\right|\mathrm{d}s\right\|_T$$
$$\leq M_1\left\|\int_{s=0}^{s=t}e^{-r(t-s)}r\left|\varphi_2(s,r)-\varphi_1(s,r)\right|\mathrm{d}s\right\|_T$$
$$\leq M_1\left\|\int_{s=0}^{s=t}e^{-r(t-s)}r\left\|\varphi_2(s,\cdot)-\varphi_1(s,\cdot)\right\|\mathrm{d}s\right\|_T \quad (91)$$
$$\leq M_1\left\|\int_{s=0}^{s=t}e^{-r(t-s)}r\left\|\varphi_2-\varphi_1\right\|_T\mathrm{d}s\right\|_T$$
$$\leq M_1\left\|\varphi_2-\varphi_1\right\|_T\sup_{0\leq t\leq T}\sup_{r>0}\int_{s=0}^{s=t}e^{-r(t-s)}r\mathrm{d}s$$
$$\leq M_1\left\|\varphi_2-\varphi_1\right\|_T$$

and

$$\left\|\int_{s=0}^{s=t}e^{-r(t-s)}r\int_{u=0}^{u=+\infty}\int_{z=0}^{z=+\infty}\left(e^{-\hat{\varphi}_2(s,u)z}-e^{-\hat{\varphi}_1(s,u)z}\right)v(\mathrm{d}z)\rho(\mathrm{d}u)\mathrm{d}s\right\|_T$$
$$\leq M_1\left\|\int_{s=0}^{s=t}e^{-r(t-s)}r\int_{u=0}^{u=+\infty}\left|\hat{\varphi}_2(s,u)-\hat{\varphi}_1(s,u)\right|\rho(\mathrm{d}u)\mathrm{d}s\right\|_T$$
$$\leq M_1\left\|\int_{s=0}^{s=t}e^{-r(t-s)}r\int_{u=0}^{u=+\infty}\left|\varphi_2(s,u)-\varphi_1(s,u)\right|\rho(\mathrm{d}u)\mathrm{d}s\right\|_T$$
$$\leq M_1\left\|\int_{s=0}^{s=t}e^{-r(t-s)}r\int_{u=0}^{u=+\infty}\left\|\varphi_2(s,\cdot)-\varphi_1(s,\cdot)\right\|\rho(\mathrm{d}u)\mathrm{d}s\right\|_T \quad (92)$$
$$= M_1\left\|\int_{s=0}^{s=t}e^{-r(t-s)}r\left\|\varphi_2(s,\cdot)-\varphi_1(s,\cdot)\right\|\mathrm{d}s\right\|_T$$
$$\leq M_1\left\|\int_{s=0}^{s=t}e^{-r(t-s)}r\left\|\varphi_2-\varphi_1\right\|_T\mathrm{d}s\right\|_T$$
$$= M_1\left\|\varphi_2-\varphi_1\right\|_T\left\|\int_{s=0}^{s=t}e^{-r(t-s)}r\mathrm{d}s\right\|_T$$
$$\leq M_1\left\|\varphi_2-\varphi_1\right\|_T$$

We therefore have a strict contraction property as well due to $0<M_1<1$ and

$$\left\|\hat{\mathbb{G}}(\varphi_1)-\hat{\mathbb{G}}(\varphi_2)\right\|_T \leq wM_1\left\|\varphi_2-\varphi_1\right\|_T+(1-w)M_1\left\|\varphi_2-\varphi_1\right\|_T=M_1\left\|\varphi_2-\varphi_1\right\|_T. \quad (93)$$

By Banach fixed point theorem (e.g., "(1.1) Theorem" in Granas and Dugundji [4]), the auxiliary equation (87) admits a unique solution. Moreover, this solution is continuous in time in any time intervals $[0,T]$ ($T>0$) because the coefficients $\theta+M_1\bar{\varphi}(\theta)$ in (88) and $M_1$ in (93) do not depend on $T$. The solution therefore exists globally. The unique solution to the auxiliary equation (87) is denoted by $\tilde{\varphi}$.

For the nonnegativity of $\tilde{\varphi}$, we have

$$\tilde{\varphi}(t,r) = \hat{\mathbb{G}}(\tilde{\varphi})(t,r) \geq 0, \ t \geq 0, \ r > 0 \tag{94}$$

since $\theta \geq 0$ and $1 - e^{-x} \geq 0$ ($x \geq 0$). Then, for any $T > 0$, the upper bound $\|\tilde{\varphi}\|_T$ of $\tilde{\varphi}$ satisfies (as for (88))

$$\begin{aligned}
\|\tilde{\varphi}\|_T &= \left\|\hat{\mathbb{G}}(\tilde{\varphi})\right\|_T \\
&= \left\|\mathbb{G}(\hat{\tilde{\varphi}})\right\|_T \\
&= \left\|\theta e^{-rt} + \int_{s=0}^{s=t} e^{-r(t-s)} r \left\{ w \int_{z=0}^{z=+\infty} \left(1 - e^{-\hat{\tilde{\varphi}}(s,r)z}\right) v(\mathrm{d}z) + (1-w) \int_{u=0}^{u=+\infty} \int_{z=0}^{z=+\infty} \left(1 - e^{-\hat{\tilde{\varphi}}(s,u)z}\right) v(\mathrm{d}z) \rho(\mathrm{d}u) \right\} \mathrm{d}s \right\|_T, \\
&\leq \theta + w \left\|\int_{s=0}^{s=t} e^{-r(t-s)} r \int_{z=0}^{z=+\infty} \left(\tilde{\varphi}(s,r) z\right) v(\mathrm{d}z) \mathrm{d}s \right\|_T \\
&\quad + (1-w) \left\|\int_{s=0}^{s=t} e^{-r(t-s)} r \int_{u=0}^{u=+\infty} \int_{z=0}^{z=+\infty} \left(\tilde{\varphi}(s,u) z\right) v(\mathrm{d}z) \rho(\mathrm{d}u) \mathrm{d}s \right\|_T \\
&\leq \theta + M_1 \|\tilde{\varphi}\|_T
\end{aligned} \tag{95}$$

yielding

$$\|\tilde{\varphi}\|_T \leq \theta(1 - M_1)^{-1} = \bar{\varphi}(q). \tag{96}$$

The upper and lower bounds of $\tilde{\varphi}$ show that it is a solution to the generalized Riccati equation (83) as well because $\tilde{\varphi} = \hat{\tilde{\varphi}}$.

Finally, we show that the generalized Riccati equation (83) admits a nonnegative unique solution $\varphi \in C_{b,T}$ with $\|\varphi\|_T \leq \bar{\varphi}(q)$ for any $T > 0$, which is $\tilde{\varphi}$. For any $T > 0$, if there exist two such solutions, then a calculation analogous to (90)-(93) yields

$$\|\varphi_2 - \varphi_1\|_T \leq M_1 \|\varphi_2 - \varphi_1\|_T \Rightarrow (1 - M_1) \|\varphi_2 - \varphi_1\|_T \leq 0 \Rightarrow \|\varphi_2 - \varphi_1\|_T = 0, \tag{97}$$

showing the uniqueness $\varphi_2 = \varphi_1$. Here, we used $M_1 \in (0,1)$.

We end this letter with a few remarks:

**Remark 1** In hindsight, the solution to the generalized Riccati equation (83) is integrable with respect to $\rho$ (because the solution is bounded continuous) as described in Proposition 3 [1].

**Remark 2** The extra parameter $\lambda$ to control the norm, which was originally used in Proof of Proposition 3 [1], is unnecessary in the presented proof strategy.

**Remark 3** If $\rho$ is a probability measure of discrete random variables (e.g., empirical measure), the supremum norm should be evaluated at each corresponding discrete point in $(0, +\infty)$.